\documentclass[reqno,10pt,centertags,draft]{amsart}
\usepackage{amsmath,amsthm,amscd,amssymb,latexsym,upref}
\date{\today}

\newcommand{\bbN}{{\mathbb{N}}}
\newcommand{\bbR}{{\mathbb{R}}}

\newcommand{\bbZ}{{\mathbb{Z}}}
\newcommand{\bbC}{{\mathbb{C}}}

\newcommand{\calD}{{\mathcal D}}

\newcommand{\calH}{{\mathcal H}}

\newcommand{\calM}{{\mathcal M}}

\newcommand{\dott}{\,\cdot\,}
\newcommand{\no}{\notag}
\newcommand{\lb}{\label}
\newcommand{\f}{\frac}

\newcommand{\ol}{\overline}

\newcommand{\wti}{\widetilde}

\newcommand{\Oh}{O}

\newcommand{\loc}{\text{\rm{loc}}}

\newcommand{\rank}{\text{\rm{rank}}}

\newcommand{\ess}{\text{\rm{ess}}}

\newcommand{\bi}{\bibitem}

\newcommand{\beq}{\begin{equation}}
\newcommand{\eeq}{\end{equation}}
\newcommand{\ba}{\begin{align}}
\newcommand{\ea}{\end{align}}

\renewcommand{\Re}{\text{\rm
Re}}
\renewcommand{\Im}{\text{\rm Im}}
\renewcommand{\ln}{\text{\rm ln}}

\renewcommand{\le}{\leqslant}

\DeclareMathOperator{\diag}{diag}

\allowdisplaybreaks
\numberwithin{equation}{section}

\newtheorem{theorem}{Theorem}[section]

\newtheorem{lemma}[theorem]{Lemma}
\newtheorem{corollary}[theorem]{Corollary}
\newtheorem{hypothesis}[theorem]{Hypothesis}
\theoremstyle{definition}
\newtheorem{definition}[theorem]{Definition}
\newtheorem{remark}[theorem]{Remark}

\begin{document}

\title[Trace Formulas and Borg-Type Theorems]{Trace Formulas and
Borg-Type Theorems \\ for Matrix-Valued Jacobi and Dirac \\ Finite
Difference Operators}
\author[S.\ Clark, F.\ Gesztesy, and W.\ Renger]{Steve Clark, Fritz
Gesztesy, and Walter Renger}
\address{Department of Mathematics and Statistics, University
of Missouri-Rolla, Rolla, MO 65409, USA}
\email{sclark@umr.edu}
\urladdr{http://www.umr.edu/\~{ }clark}
\address{Department of Mathematics,
University of Missouri, Columbia, MO 65211, USA}
\email{fritz@math.missouri.edu} 
\urladdr{http://www.math.missouri.edu/people/fgesztesy.html}
\address{Dr. Johannes Heidenhain GmbH, 83301 Traunreut, Germany}
\email{walter\_renger@web.de}
\thanks{Based upon work supported by the US National Science
Foundation under Grants No.\ DMS-0405526 and DMS-0405528.}
\begin{abstract}
Borg-type uniqueness theorems for matrix-valued Jacobi operators
$H$ and supersymmetric Dirac difference operators $D$ are proved. More
precisely, assuming reflectionless matrix coefficients $A, B$ in the 
self-adjoint Jacobi operator $H=AS^++A^-S^-+B$ (with $S^\pm$ the
right/left shift operators on the lattice $\bbZ$) and the 
spectrum of $H$ to be a compact interval $[E_-,E_+]$, $E_-<E_+$, we
prove that $A$ and $B$ are certain multiples of the identity matrix. An
analogous result which, however,  displays a certain novel
nonuniqueness feature, is proved for supersymmetric self-adjoint Dirac
difference operators $D$ with spectrum  given by
$\big[-E_+^{1/2},-E_-^{1/2}\big]\cup\big[E_-^{1/2},E_+^{1/2}\big]$, 
$0 \leq E_- < E_+$.

Our approach is based on trace formulas and matrix-valued
(exponential) Herglotz representation theorems. As a by-product 
of our techniques we obtain the extension of Flaschka's Borg-type 
result for periodic scalar Jacobi operators to the class of
reflectionless matrix-valued Jacobi operators.
\end{abstract}

\maketitle

\section{Introduction}\label{s1}

As discussed in detail in \cite{GKM02}, while various aspects of
inverse spectral theory for scalar  Schr\"odinger, Jacobi, and
Dirac-type operators, and more generally, for $2\times 2$ Hamiltonian
systems, are  well-understood by now, the corresponding theory for
such  operators and Hamiltonian systems with $m\times m$, $m\in\bbN$, 
matrix-valued coefficients is still largely a wide open field. A
particular  inverse spectral theory aspect we have in mind is that of
determining isospectral sets (manifolds) of such systems. In this
context it may,  perhaps, come as a surprise that even determining the
isospectral set  of Hamiltonian systems with matrix-valued periodic
coefficients  is an open problem. The present paper is a modest
attempt  toward a closer investigation of inverse spectral problems in
connection with Borg-type uniqueness theorems for matrix-valued Jacobi
operators and Dirac-type finite difference systems. It should be
mentioned that these problems are not just of interest in a spectral
theoretic context, but due to their implications for other areas such
as completely integrable systems (e.g., the nonabelian Toda and
Kac--van Moerbeke hierarchies), are of interest to a larger audience. 

Before we describe the content of this paper in more detail, we
briefly comment on background literature for matrix-valued Jacobi and
Dirac-type difference operators. Spectral and Weyl--Titchmarsh theory 
for Jacobi operators can be found in \cite[Sect.\ VII.2]{Be68},
\cite{Fu76}, \cite{RS02}, \cite[Ch.\ 10]{Sa97} and the literature
therein. The case of Dirac finite difference operators was discussed
in detail in \cite{CG03}. Deficiency indices of matrix-valued Jacobi
operators are studied in \cite{KM98}--\cite{KM01}. Inverse spectral and
scattering theory for matrix-valued finite difference systems and its
intimate connection to matrix-valued orthogonal polynomials and the
moment problem are treated in \cite{AG94}, \cite{AN84}, \cite[Sect.\
VII.2]{Be68}, \cite{DL96}--\cite{DV95}, \cite{Ge82}, \cite{Lo99},
\cite{MBO92}, \cite{Os00}, \cite{Os02}, \cite{RS02}, \cite[Ch.\
8]{Sa97}, \cite{Se80}. A number of uniqueness theorems for
matrix-valued Jacobi operators were proved in \cite{GKM02}. Finally,
connections with  nonabelian completely integrable systems are
discussed in \cite{BGS86}, \cite{BG90}, \cite{Os97}, \cite{Sa03},
\cite[Chs.\ 9, 10]{Sa97}.

Let $\bbC(\bbZ)^{r\times s}$ be the space of sequences of complex 
$r\times s$ matrices, $r,s\in\bbN$. We also denote 
$\bbN_0=\bbN\cup\{0\}$ in the following.

The matrix-valued Jacobi operators $H$ in $\ell^2(\bbZ)^m$ discussed
in this paper are of the form 
\begin{equation} 
H = A S^+ + A^- S^-  + B, \quad \calD(H)=\ell^2(\bbZ)^m. \label{1.1} 
\end{equation}
Here $S^\pm$ denote the shift operators acting upon
$\bbC(\bbZ)^{r\times s}$, $r,s\in\bbN$, as
\begin{equation}
S^{\pm}f(\,\cdot\,)=f^{\pm}(\,\cdot\,)=f (\, \cdot\,\pm 1), \quad 
f\in\bbC(\bbZ)^{r\times s}, \lb{1.2}
\end{equation}
and $A=\{A(k)\}_{k\in\bbZ}\in \bbC(\bbZ)^{m\times m}$, 
$A(k)>0$ for all $k\in\bbZ$, $B=\{B(k)\}_{k\in\bbZ}\in
\bbC(\bbZ)^{m\times m}$, $B(k)=B(k)^*$ for all $k\in\bbZ$ and
$m\in\bbN$. 

The Dirac-type finite difference operators $D$ in 
$\ell^2(\bbZ)^m \oplus \ell^2(\bbZ)^m$ studied in this paper are of
the form 
\begin{equation} \label{1.3}
D=S_{\rho} + X= \begin{pmatrix}0 & \rho
S^+ + \chi^*\\ \rho^- S^- + \chi & 0 \end{pmatrix}, \quad 
\calD(D)=\ell^2(\bbZ)^m \oplus \ell^2(\bbZ)^m, 
\end{equation}
where $\rho=\{\rho(k)\}_{k\in\bbZ}\subset\bbC(\bbZ)^{m\times m}$ 
and $S_\rho$ and X are of the
block form
\begin{equation}
S_{\rho}= \begin{pmatrix}0 & \rho
S^+  \\ \rho^- S^- & 0 \end{pmatrix}, \quad 
X=\begin{pmatrix} 0 & \chi^* \\ \chi & 0 \end{pmatrix} \lb{1.4}
\end{equation}
with $\rho(k), \chi(k)\in\bbC^{m\times m}$ invertible for all
$k\in\bbZ$.  Moreover, following \cite{CG03}, we may assume
without loss of generality that for all $k\in\bbZ$, $\rho(k)$ is a
positive definite diagonal $m\times m$ matrix (cf.\ Remark\ \ref{r4.12a}
for details)). 

Additional assumptions on the coefficients $A(k)$, $B(k)$, $\rho(k)$,
and $\chi(k)$, $k\in\bbZ$, will be formulated in Sections \ref{s2} and
\ref{s5}, respectively. 

The difference operators $H$ and $D$ represent the natural
matrix-valued generalizations of Lax operators arising in connection
with Kac-van Moerbeke and Toda lattices  (cf.\ \cite{BGHT98},
\cite{GH05} and the references therein) and hence lead to nonabelian
Toda and Kac--van Moerbeke hierarchies of completely integrable
nonlinear evolution equations.  

Next we briefly describe the history of Borg-type theorems relevant to
this paper. In 1946, Borg \cite{Bo46} proved, among a variety of other 
inverse spectral theorems, the following result. (We denote by 
$\sigma(\cdot)$ and $\sigma_{\ess}(\cdot)$ the spectrum and essential
spectrum of a densely defined closed linear operator.)

\begin{theorem}[\cite{Bo46}] \lb{t1.1}
Let $q\in L^1_{\loc} (\bbR)$ be real-valued and periodic. Let
$h=-\f{d^2}{dx^2}+q$ be the associated self-adjoint 
Schr\"odinger operator in $L^2(\bbR)$ and suppose that
$\sigma(h)=[e_0,\infty)$ for some $e_0\in\bbR$.  Then $q$ is of the
form,
\begin{equation}
q(x)=e_0 \, \text{ for a.e.\ $x\in\bbR$}. \lb{1.5}
\end{equation}
\end{theorem}

\begin{remark}\lb{r1.1}
Traditionally, uniqueness results such as Theorem\ \ref{t1.1} are
called Borg-type theorems. This terminology, although generally
accepted, is a bit unfortunate as the same term is also used for
other theorems Borg proved in his celebrated 1946 paper
\cite{Bo46}. Indeed, inverse spectral results on finite intervals in
which the
potential coefficient(s) are recovered from several spectra were also
pioneered by Borg in \cite{Bo46} and theorems of this kind are now also
described as Borg-type theorems in the literature, see, e.g.,
\cite{Ma94}--\cite{Ma99a}.
\end{remark}

A closer examination of the proof of Theorem\ \ref{t1.1} in
\cite{CGHL00} shows that periodicity of $q$ is not the point for the
uniqueness result \eqref{1.5}. The key ingredient (besides
$\sigma(h)=[e_0,\infty)$ and $q$ real-valued) is the fact that 
\begin{equation}
\text{for all $x\in \bbR$, } \, \xi(\lambda,x)=1/2 \, \text{ for a.e.\ 
$\lambda\in\sigma_{\ess}(h)$.} \lb{1.6}
\end{equation}
Here $\xi(\lambda,x)$, the argument of the boundary value
$g(\lambda+i0,x)$ of the diagonal Green's function of $h$ on the real
axis (where $g(z,x)=(h-zI)^{-1}(x,x))$, $z\in\bbC\backslash\sigma(h)$),
is defined by 
\begin{equation}
\xi(\lambda,x)=\pi^{-1}\lim_{\varepsilon\downarrow 0} 
\Im(\ln(g(\lambda+i\varepsilon,x)) \, \text{  
for a.e.\ $\lambda\in\bbR$}. \lb{1.7}
\end{equation}

Real-valued periodic potentials are known to satisfy \eqref{1.6}, but
so do certain classes of real-valued quasi-periodic and
almost-periodic potentials $q$. In particular, the class of real-valued
algebro-geometric finite-gap KdV potentials $q$ (a subclass of the set
of real-valued quasi-periodic  potentials) is a prime example
satisfying \eqref{1.6} without necessarily being periodic. 
Traditionally, potentials $q$ satisfying \eqref{1.6} are called
\textit{reflectionless} (see \cite{CG02}, \cite{CGHL00} and the
references therein).

The extension of Borg's Theorem\ \ref{t1.1} to periodic matrix-valued
Schr\"odinger operators was proved by D\'epres \cite{De95}. A new
strategy of the proof based on exponential Herglotz representations
and a trace formula for such potentials, as well as the extension to
reflectionless matrix-valued potentials, was obtained in
\cite{CGHL00}.

The direct analog of Borg's Theorem\ \ref{t1.1} for scalar (i.e.,
$m=1$) periodic Jacobi operators was proved by Flaschka \cite{Fl75} in
1975.

\begin{theorem}  [\cite{Fl75}] \lb{t1.2}
Suppose $a$ and $b$ are periodic real-valued sequences in
$\ell^\infty(\bbZ)$ with the same period and $a(k)>0$, $k\in\bbZ$. Let
$h=aS^+ +a^-S^- +b$ be the associated self-adjoint  Jacobi operator
in $\ell^2(\bbZ)$ and suppose that $\sigma(h)=[E_-,E_+]$ for some
$E_-<E_+$.  Then $a$ and $b$ are of the form,
\begin{equation}
a(k)=(E_+-E_-)/4 , \quad  b(k)=(E_-+E_+)/2, \quad k\in\bbZ.  \lb{1.8}
\end{equation}
\end{theorem}

The extension of Theorem\ \ref{t1.2} to reflectionless scalar Jacobi
operators is due to Teschl \cite[Corollary\ 6.3]{Te98} (see also
\cite[Corollary\ 8.6]{Te00}). As one of the principal results in this
paper we will extend Theorem\ \ref{t1.2} to matrix-valued
reflectionless Jacobi operators $H$ of the type \eqref{1.1} in Section\
\ref{s4}. This extension (and the special case of periodic
matrix-valued Jacobi operators) is new. In addition, we will prove a
Borg-type theorem for the supersymmetric Dirac difference operator
$D$ in \eqref{1.3} which is new even in the simplest case $m=1$. The
latter displays an interesting nonuniqueness feature which has not
previously been encountered with Borg-type theorems.

In Section\ \ref{s2} we review the basic Weyl--Titchmarsh theory and
the corresponding Green's matrices for matrix-valued Jacobi operators
on $\bbZ$ and on a half-lattice. Section\ \ref{s3} is devoted to
asymptotic expansions of Weyl--Titchmarsh and Green's matrices as the
(complex) spectral parameter tends to infinity. 
Section\ \ref{s4} contains the derivation of a trace formula for
Jacobi operators and one of our principal new results, the
proof of a Borg-type theorem for matrix-valued Jacobi operators. 
Finally, Section\ \ref{s5} presents a Borg-type  theorem for
supersymmetric Dirac-type difference operators.

\section{Weyl--Titchmarsh and Green's Matrices for Matrix-Valued \\
Jacobi Operators} \label{s2}

In this section we consider Weyl--Titchmarsh and Green's matrices 
for matrix-valued Jacobi operators on $\bbZ$ and on a half-lattice. 

We closely follow the treatment of matrix-valued Jacobi operators
in \cite{GKM02}. As the basic hypothesis in this section we adopt the
following set of asumptions.

\begin{hypothesis} \lb{h4.1}
Let $m\in\bbN$ and consider the sequences of 
self-adjoint $m\times m$ matrices
\begin{align}
&A=\{A(k)\}_{k\in\bbZ}\in \bbC(\bbZ)^{m\times m}, \quad 
A(k)>0, \,\, k\in\bbZ, \no \\
&B=\{B(k)\}_{k\in\bbZ}\in \bbC(\bbZ)^{m\times m}, \quad B(k)=B(k)^*,
\,\, k\in\bbZ. \lb{4.1}
\end{align}
Moreover, assume that $A(k)$ and $B(k)$ are uniformly bounded
with  respect to $k\in\bbZ$, that is, there exists a $C>0$,
such that 
\begin{equation}
\|A(k)\|_{\bbC^{m\times m}}+\|B(k)\|_{\bbC^{m\times m}}\leq C,
\quad k\in\bbZ. \lb{4.1a}
\end{equation}
\end{hypothesis}
Next, denote by $S^\pm$ the shift operators in 
$\bbC(\bbZ)^{r\times s}$, $r,s\in\bbN$, 
\begin{equation}
(S^\pm g)(k)=g^\pm (k)=g(k\pm 1), \quad 
g\in \bbC (\bbZ)^{r\times s}, 
\,\,k\in\bbZ.  \lb{4.7}
\end{equation}
Given Hypothesis\ \ref{h4.1}, the matrix-valued self-adjoint Jacobi 
operator $H$ with domain $\ell^2(\bbZ)^m$ is then defined by
\begin{equation} 
H = A S^+ + A^- S^-  + B, \quad \calD(H)=\ell^2(\bbZ)^m.  
\label{4.8b} 
\end{equation}
Because of hypothesis \eqref{4.1a}, $H$ is a bounded symmetric
operator and hence self-adjoint. In particular, the difference 
expression $A S^+ + A^- S^-  + B$ induced by \eqref{4.8b} is in the
limit point case at $\pm\infty$. We chose to adopt \eqref{4.1a} for
simplicity only. Our formalism extends to unbounded Jacobi
operators and to cases  where the difference expression associated
with \eqref{4.8b} is in  the limit circle case at $+\infty$ and/or
$-\infty$ (cf.\ \cite{CG03}).  We just note in passing that without
assuming \eqref{4.1a}, the difference expression $A S^+ + A^- S^-  + B$
is in the limit point case at $\pm\infty$ if
$\sum_{k=k_0}^{\pm\infty}\|A(k)\|_{\bbC^{m\times m}}^{-1} =\infty$ (see,
e.g., \cite[Theorem\ VII.2.9]{Be68}). Whenever we need to stress the
dependence of $H$ on $A, B$ we will write $H(A,B)$ instead of $H$. 

The Green's matrix associated with $H$ will be denoted by $G(z,k,\ell)$
in the following,
\begin{equation}
G(z,k,\ell)=(H-zI_m)^{-1}(k,\ell), \quad z\in\bbC\backslash\bbR, \,\, 
k,\ell\in\bbZ. \lb{4.8a}
\end{equation}

Next, fix a site $k_0 \in \bbZ$ and define $m\times m$  
matrix-valued solutions $\phi(z,k,k_0)$ and
$\theta(z,k,k_0)$ of the equation
\begin{equation}
A(k)\psi(z,k+1)+A(k-1)\psi(z,k-1)+(B(k)-zI_m)\psi(z,k)=0, 
\quad z\in\bbC, \, k\in\bbZ, \lb{4.2}
\end{equation}
satisfying the initial conditions
\begin{equation}\label{4.3}
\theta(z,k_0,k_0)=\phi(z,k_0+1,k_0)=I_m, \quad 
\phi(z,k_0,k_0)=\theta(z,k_0+1,k_0)=0.
\end{equation}
One then introduces $m\times m$ matrix-valued Weyl--Titchmarsh solutions
$\psi_{\pm}(z,k,k_0)$ associated with $H$ defined by 
\begin{equation}
\psi_{\pm}(z,k,k_0)= \theta (z,k,k_0) - \phi (z,k,k_0)A(k_0)^{-1}
M_\pm(z,k_0), \quad z\in\bbC\backslash\bbR, \,\, k\in\bbZ, \lb{4.4}
\end{equation}
with the properties 
\begin{align}
&A(k)\psi_{\pm}(z,k+1,k_0)+A(k-1)\psi_\pm(z,k-1,k_0)
+(B(k)-zI_m)\psi_\pm(z,k,k_0)=0, \no \\
& \hspace*{8.7cm} z\in\bbC\backslash\bbR, \,\, k\in\bbZ, 
\lb{4.4a} \\ 
&\psi_\pm (z,\cdot,k_0)\in
\ell^2((k_0,\pm\infty)\cap
\bbZ)^{m\times m}, \quad z\in\bbC\backslash\bbR, \lb{4.5} 
\end{align}
where $M_\pm(z,k_0)$ denote the half-line Weyl--Titchmarsh matrices
associated with $H$. Since by assumption $A S^+ + A^- S^- + B$ 
is in the limit point case at $\pm\infty$, $M_\pm(z,k_0)$ in
\eqref{4.4} are uniquely determined by the requirement \eqref{4.5}.
We also note that by a standard argument, 
\begin{equation}
\det(\phi(z,k,k_0))\neq 0 \text{ for all }
k\in\bbZ\backslash\{k_0\} \text{ and } z\in\bbC\backslash\bbR, 
\lb{4.4b}
\end{equation}
since otherwise one could construct a Dirichlet-type
eigenvalue $z\in\bbC\backslash\bbR$ for $H$ restricted to the finite
segment $k_0+1,\dots,k-1$ for
$k\geq k_0+1$ of $\bbZ$ (and similarly for $k\leq k_0-1$).
Thus, introducing 
\begin{equation}
M_N(z,k_0)=-\phi(z,N,k_0)^{-1}\theta(z,N,k_0), \quad
z\in\bbC\backslash\bbR, \,\, N\in\bbZ\backslash\{k_0\}, \lb{4.5a}
\end{equation}
one can then compute $M_\pm(z,k_0)$ by the limiting relation
\begin{equation}
M_\pm(z,k_0)=\lim_{N\to\pm\infty} M_N(z,k_0), \quad
z\in\bbC\backslash\bbR, \lb{4.5b}
\end{equation}
the limit being unique since $A S^+ + A^- S^-  + B$ is in
the limit point case at $\pm\infty$. Alternatively, \eqref{4.4}
yields
\begin{equation}
M_\pm(z,k_0)=-A(k_0)\psi_\pm(z,k_0+1,k_0), \quad
z\in\bbC\backslash\bbR. \lb{4.5c}
\end{equation}
More generally, recalling
\begin{equation}
\det(\psi_\pm(z,k,k_0))\neq 0 \text{ for all }
k\in\bbZ \text{ and } z\in\bbC\backslash\bbR, 
\lb{4.5d}
\end{equation}
by an argument analogous to that following \eqref{4.4b}, we now 
introduce
\begin{equation}
M_\pm(z,k)=-A(k)\psi_\pm(z,k+1,k_0)\psi_\pm(z,k,k_0)^{-1}, \quad
z\in\bbC\backslash\bbR, \,\, k\in\bbZ. \lb{4.5e}
\end{equation}
One easily verifies that $M_\pm(z,k)$ represent the
Weyl--Titchmarsh $M$-matrices associated with the reference point
$k\in\bbZ$. Moreover, one obtains the Riccati-type equation for
$M_\pm(z,k)$, 
\begin{equation}
M_\pm(z,k)+A(k-1)M_\pm(z,k-1)^{-1}A(k-1)+zI_m-B(k)=0, \quad 
z\in\bbC\backslash\bbR, \,\, k\in\bbZ \lb{4.5f}
\end{equation}
as a result of \eqref{4.4a}. For later reference we summarize the
principal results on $M_{\pm}(z,k_0)$ in the following theorem. (We
denote as usual $\Re(M)=(M+M^*)/2$, $\Im(M)=(M^*-M)/(2i)$, etc., for
square matrices $M$.)
\begin{theorem}[\cite{AD56}, \cite{Ca76}, \cite{GT00}--\cite{HS86},
\cite{KS88}] \lb{thm2.3}  Assume Hypothesis\ \ref{h4.1} and suppose
that   
$z\in\bbC\backslash\bbR$, and $k_0\in\bbZ$. Then, \\
$(i)$ $\pm M_{\pm}(z,k_0)$ is a matrix-valued Herglotz 
function of maximal rank. In particular,
\begin{gather}
\Im(\pm M_{\pm}(z,k_0)) > 0, \quad z\in\bbC_+, 
\lb{4.6A} \\
M_{\pm}(\bar z,k_0)=M_{\pm}( z,k_0)^*,\lb{4.6a} \\
\rank (M_{\pm}(z,k_0))=m, \lb{4.6b} \\
\lim_{\varepsilon\downarrow 0} M_{\pm}(\lambda
+i\varepsilon,k_0) \text{ exists for a.e.\
$\lambda\in\bbR$}. \lb{4.6c}
\end{gather}
Isolated poles of $\pm M_{\pm}(z,k_0)$ and $\mp
M_{\pm}(z,k_0)^{-1}$ are at most of first order, 
are real, and have a nonpositive residue.  \\
$(ii)$  $\pm M_{\pm}(z,k_0)$ admit the representations
\begin{align}
\pm M_{\pm}(z,k_0)&=K_\pm(k_0)+L_{\pm}(k_0)z+\int_\bbR 
d\Omega_\pm(\lambda,k_0) \,
\bigg(\f{1}{\lambda-z}-\f{\lambda}{1+\lambda^2}\bigg) 
\lb{4.6d}  \\
&=\exp\bigg[C_\pm(k_0)+\int_\bbR d\lambda \, 
\Xi_{\pm} (\lambda,k_0)
\bigg(\f{1}{\lambda-z}-\f{\lambda}{1+\lambda^2}\bigg)  
 \bigg], \lb{4.6e}
\end{align}
where
\begin{align}
K_\pm(k_0)&=K_\pm(k_0)^*, \; L_+(k_0)=0, \; L_-(k_0)=L_-(k_0)^*, 
\quad \int_\bbR \f{\Vert
d\Omega_\pm(\lambda,k_0)\Vert}{1+\lambda^2}<\infty, 
\lb{4.6f} \\
C_\pm(k_0)&=C_\pm(k_0)^*, 
\quad 0\le\Xi_\pm(\dott,k_0)\le I_m \; \rm{  a.e.}
\lb{4.6g}
\end{align}
Moreover,
\begin{align}
\Omega_\pm((\lambda,\mu],k_0)&=\lim_{\delta\downarrow
0}\lim_{\varepsilon\downarrow 0}\f1\pi
\int_{\lambda+\delta}^{\mu+\delta} d\nu \, \Im(\pm
M_\pm(\nu+i\varepsilon,k_0)), \lb{4.6h} \\
\Xi_\pm(\lambda,k_0)&=
\lim_{\varepsilon\downarrow 0}\f1\pi\Im(\ln(\pm
M_\pm(\lambda+i\varepsilon,k_0)) \text{ for a.e.\ 
$\lambda\in\bbR$}.\lb{4.6i}
\end{align}
\end{theorem}

Next, we define the self-adjoint half-line Jacobi operators 
$H_{\pm,k_0}$ on $\ell^2([k_0,\pm\infty)\cap\bbZ)^m$ by 
\begin{equation}\lb{4.9}
H_{\pm, k_0}=P_{\pm, k_0}HP_{\pm,
k_0}\big|_{\ell^2((k_0,\pm\infty)\cap\bbZ)^m}, \quad k_0\in \bbZ,
\end{equation}
where $P_{\pm, k_0}$ are the orthogonal
projections onto the subspaces
$\ell^2([k_0,\pm \infty)\cap \bbZ)^m$. In addition, Dirichlet boundary
conditions at $k_0\mp 1$ are associated with $H_{\pm,k_0}$, 
\begin{align} 
(H_{+, k_0}f)(k_0) &= A(k_0)f^+(k_0)+B(k_0)f(k_0), \no \\
(H_{-, k_0}f)(k_0) &= A^-(k_0)f^-(k_0)+B(k_0)f(k_0), \no \\
(H_{\pm, k_0}f)(k) &= A(k)f^+(k)+A^-(k)f^-(k)+B(k)f(k), \quad 
k\gtreqless k_0\pm 1, \label{4.8}  \\
f\in \calD(H_{\pm, k_0})&=\ell^2([k_0,\pm\infty)\cap\bbZ)^m \no
\end{align}
(i.e., formally, $f(k_0\mp 1)=0$). We also introduce  
$m$-functions $m_\pm(z,k_0)$ associated with $H_{\pm, k_0}$ by
\begin{align}
m_\pm(z,k_0)&=Q_{k_0}(H_{\pm, k_0}-zI_m)^{-1}Q_{ k_0}, \lb{4.10} \\
&=G_{\pm,k_0}(z,k_0,k_0), \quad z\in\bbC\backslash\bbR,
\,\, k_0\in \bbZ. \lb{4.10a}
\end{align}
Here $Q_{ k_0}$ are orthogonal projections
onto the $m$-dimensional subspaces $\ell^2(\{k_0 \})^m$,
$k_0\in\bbZ$  and 
\begin{equation}
G_{\pm,k_0}(z,k,\ell)=(H_{\pm,k_0}-zI_m)^{-1}(k,\ell), \quad
z\in\bbC\backslash\bbR, \,\,  k,\ell\in\bbZ\cap [k_0,\pm\infty)
\lb{4.10b}
\end{equation}
represent the Green's matrices of $H_{\pm,k_0}$.

In order to find the connection between $m_\pm(z,k_0)$ and
$M_\pm(z,k_0)$ we briefly discuss the Green's matrices
$G_{\pm,k_0}(z,k,\ell)$ and $G(z,k,\ell)$ associated
with $H_{\pm, k_0}$ and $H$ next. 

First we recall the definition of the Wronskian $W(f,g)(k)$ of two
sequences of matrices $f(\cdot), g(\cdot)\in \bbC(\bbZ)^{m\times m} $
given by 
\begin{equation}\lb{4.6}
W(f,g)(k)=f(k)A(k)g(k+1)-f(k+1)A(k)g(k), \quad k\in\bbZ.
\end{equation}
We note that for any two matrix-valued solutions 
$\varphi(z,\cdot)$ and  $\psi(z,\cdot)$ of \eqref{4.2}
the Wronskian $W(\varphi(\overline{ z},\cdot)^*,\psi( z,\cdot))(k)$
is independent of $k\in\bbZ$.

In complete analogy to the scalar Jacobi case (i.e., $m=1$) one
verifies,
\begin{align}
&G_{+, k_0}(z,k,\ell)=\begin{cases} -\psi_+(z,k,k_0-1) 
A(k_0-1)^{-1}\phi(\ol z, \ell,k_0-1)^*, & \ell\leq k, \\
-\phi(z,k,k_0-1) 
A(k_0-1)^{-1}\psi_+(\ol z, \ell,k_0-1)^*, & \ell\geq k, 
\end{cases} \lb{4.10c} \\
& \hspace*{6.55cm} z\in\bbC\backslash\bbR, \,\,  k,\ell\in\bbZ\cap
[k_0,\infty). \no \\
&G_{-, k_0}(z,k,\ell)=\begin{cases} \phi(z,k,k_0+1) 
A(k_0+1)^{-1}\psi_-(\ol z, \ell,k_0+1)^*, & \ell\leq k, \\
\psi_-(z,k,k_0+1) 
A(k_0+1)^{-1}\phi(\ol z, \ell,k_0+1)^*, & \ell\geq k, 
\end{cases} \lb{4.10c'} \\
& \hspace*{6cm} z\in\bbC\backslash\bbR, \,\,  k,\ell\in\bbZ\cap
(-\infty,k_0]. \no
\end{align}
Similarly, using the fact that 
\begin{align}
&\psi_+(z,k,k_0)[M_-(z,k_0)
-M_+(z,k_0)]^{-1}\psi_- (\overline {z},k,k_0)^* \no \\
&= \psi_-(z,k,k_0)[M_-(z,k_0)
-M_+(z,k_0)]^{-1}\psi_+ (\overline{z},k,k_0)^* \no \\
&=[M_-(z,k)-M_+(z,k)]^{-1}, \quad k\in \bbZ, \lb{4.16}
\end{align}
and 
\begin{align}
&A(k)\psi_+(z,k+1,k_0)[M_-(z,k_0)
-M_+(z,k_0)]^{-1}\psi_- (\overline z,k,k_0)^*  \lb{4.17} \\
&- A(k)\psi_-(z,k+1,k_0)[M_-(z,k_0)-M_+(z,k_0)]^{-1}
\psi_+ (\overline z,k,k_0)^*=I_m, \quad k\in \bbZ, \no
\end{align}
one verifies 
\begin{align}
&G(z,k,\ell) = \begin{cases} \psi_+(z,k,k_0)[M_-(z,k_0)
-M_+(z,k_0)]^{-1}\psi_- (\overline z,\ell,k_0)^*,  & \ell \leq k,
 \\
\psi_-(z,k,k_0)[M_-(z,k_0)
-M_+(z,k_0)]^{-1}\psi_+ (\overline z,\ell,k_0)^*,  & \ell \geq k,
\end{cases} \label{4.15} \\
& \hspace*{8.3cm} z\in\bbC\backslash\bbR, \,\,  k,\ell\in\bbZ. \no
\end{align}

Using \eqref{4.3}, \eqref{4.4a}, \eqref{4.5f}, and \eqref{4.10a}, 
one infers that the Weyl--Titchmarsh matrices $M_\pm(z,k)$
introduced by \eqref{4.5b} (resp.~\eqref{4.5c}) and the 
$m$-functions $m_\pm(z,k)$ defined in \eqref{4.10} are 
related by
\begin{align} 
M_+(z,k)&=-m_+(z, k)^{-1}-zI_m+B(k), \quad z\in\bbC\backslash\bbR, \;
k\in\bbZ 
\lb{4.11}
\intertext{and}
M_-(z,k)&=m_-(z, k)^{-1}, \quad z\in\bbC\backslash\bbR, \;
k\in\bbZ .  \lb{4.12}
\end{align}

In analogy to \eqref{4.5f}, $m_\pm(z,k)$ also satisfy Riccati-type
equations of the form
\begin{align}
&A(k-1)m_+(z,k)A(k-1)m_+(z,k-1) +(zI_m-B(k-1))m_+(z,k-1)+I_m=0, \no \\
& \hspace*{8.5cm} z\in\bbC\backslash\bbR, \; k\in\bbZ \lb{4.12a}
\intertext{and}
&A(k-1)m_-(z,k-1)A(k-1)m_-(z,k) +(zI_m-B(k))m_-(z,k)+I_m=0, \no \\
& \hspace*{8.5cm} z\in\bbC\backslash\bbR, \; k\in\bbZ.  \lb{4.12b}
\end{align}

Next, we introduce the $2m\times 2m$ Weyl--Titchmarsh matrix 
$M(z,k)$ associated with the Jacobi operator $H$ in 
$\ell^2(\bbZ)^m$ by
\begin{equation}
M(z,k)=\big(M_{j,j^\prime}(z,k)\big)_{j,j^\prime=1,2}\, , 
\quad z\in\bbC\backslash\bbR, \,\, k\in\bbZ, \lb{4.21}
\end{equation}
where
\begin{align}
M_{1,1}(z,k)&=[M_{-}(z,k)-M_{+}(z,k)]^{-1}, \lb{4.22} \\
M_{1,2}(z,k)&=2^{-1}[M_{-}(z,k)-M_{+}(z,k)]^{-1}
[M_{-}(z,k)+M_{+}(z,k)], \lb{4.23} \\
M_{2,1}(z,k)&=2^{-1}[M_{-}(z,k)+M_{+}(z,k)]
[M_{-}(z,k)-M_{+}(z,k)]^{-1}, \lb{4.23a} \\
M_{2,2}(z,k)&=M_{\pm}(z,k)
[M_{-}(z,k)-M_{+}(z,k)]^{-1}M_{\mp}(z,k). \lb{4.23b} 
\end{align}

One verifies that $M(z,k)$ is a $2m\times 2m$ Herglotz matrix with the
following properties: 

\begin{theorem} [\cite{AD56}, \cite{Ca76}, \cite{GT00}--\cite{HS86}]
\lb{thm2.4}  Assume Hypothesis\ \ref{h4.1}, $z\in\bbC \backslash
\bbR$, and
$k_0\in\bbZ$.  Then, $M(z,k_0)$ is a matrix-valued Herglotz function of 
rank $2m$ with representations
\begin{align}
M(z,k_0)&=K(k_0)+\int_\bbR d\Omega(\lambda,k_0)\,
\bigg(\f1{\lambda-z}-\f{\lambda}{1+\lambda^2}\bigg)
\lb{4.23c} \\
&=\exp\bigg[C(k_0)+\int_\bbR d\lambda\, 
\Upsilon (\lambda,k_0)
\bigg(\f1{\lambda-z}-\f{\lambda}{1+\lambda^2} \bigg) \bigg], \lb{4.23d}
\end{align}
where
\begin{align}
K(k_0)&=K(k_0)^*, \quad \int_\bbR\f{\Vert 
d\Omega(\lambda,k_0)
\Vert}{1+\lambda^2}<\infty,
\lb{4.23e}\\
C(k_0)&=C(k_0)^*, \quad 0\le \Upsilon (\dott,k_0)\le I_{2m} 
\text{ a.e.} \lb{4.23f}
\end{align}
Moreover,
\begin{align}
\Omega((\lambda,\mu],k_0)&=\lim_{\delta\downarrow
0}\lim_{\varepsilon\downarrow 0}\f1\pi
\int_{\lambda+\delta}^{\mu+\delta} d\nu \, \Im(
M(\nu+i\varepsilon,k_0)), \lb{4.23g} \\
\Upsilon(\lambda,k_0)&=\lim_{\varepsilon\downarrow
0}\f1\pi\Im(\ln(M(\lambda+i\varepsilon,k_0)))
\text{ for a.e.\ $\lambda\in\bbR$}.\lb{4.23h}
\end{align}
\end{theorem}

\begin{remark}
We note that the Weyl--Titchmarsh matrix 
$M(z,k)$ is related to the Green's matrix associated with the 
Jacobi operator $H$ by
\begin{equation}
M(z,k)=\begin{pmatrix} I_m &0  \\
0 &-A(k) \end{pmatrix} \calM(z,k)
\begin{pmatrix}I_m &0  \\
0 &-A(k) \end{pmatrix}
+\frac{1}{2}\begin{pmatrix} 0 &I_m  \\
I_m &0 \end{pmatrix},  \lb{4.23i}
\end{equation}
where
\begin{equation}
\calM(z,k)=
\begin{pmatrix}G(z, k,k) &G(z,k,k+1)  \\
 G(z,k+1, k) & G(z, k+1, k+1) \end{pmatrix}.  \lb{4.23j}
\end{equation} 
\end{remark}

With $M(z,k)$ defined in \eqref{4.21}, the following 
uniqueness theorem holds. 

\begin{theorem} \lb{t4.2}
Assume Hypothesis\ \ref{h4.1} and let $k_0\in\bbZ$. Then
the $2m\times 2m$ Weyl--Titchmarsh matrix 
$M(z,k_0)$ for all $z\in\bbC_+$ uniquely determines the 
Jacobi operator $H$ and hence $A=\{A(k)\}_{k\in\bbZ}$ and 
$B=\{B(k)\}_{k\in\bbZ}$.
\end{theorem}

Perhaps the simplest way to prove Theorem\ \ref{t4.2} is to reduce
it  to knowledge of $M_\pm(z,k_0)$, and hence by
\eqref{4.11} and \eqref{4.12} to that of $m_\pm(z,k_0)$ 
 for all $z\in\bbC_+$. (We note that the knowledge of $B(k_0)$, which 
is required according to \eqref{4.11}, can be determined from the
asymptotics of $M_-(z,k_0)$ in \eqref{3.8}.) Using the standard
construction of orthogonal matrix-valued polynomials with respect to the
normalized measure in the Herglotz representation of
$m_\pm(z,k_0)$,
\begin{equation}
 m_\pm(z,k_0)=\int_\bbR d\nu_\pm(\lambda,k_0) \,
(\lambda-z)^{-1}, \quad z\in\bbC_+, \quad
\int_\bbR d\nu_\pm(\lambda,k_0) =I_m, \lb{4.24}
\end{equation}
allows one
to reconstruct $A(k)$, $B(k)$, $k\in [k_0,\pm\infty)\cap\bbZ$
from the measures $d\nu_\pm(\lambda,k_0)$ (cf., e.g., 
\cite[Section~VII.2.8]{Be68}).
More precisely,
\begin{align}
&A(k_0\pm k)=\int_\bbR\lambda P_{\pm,k}(\lambda, k_0)
 d\nu_\pm(\lambda,k_0) \,
P_{\pm,k+1}(\lambda, k_0)^*,\no \\
& B(k_0\pm k)=\int_\bbR\lambda P_{\pm,k}(\lambda,k_0)
 d\nu_\pm(\lambda,k_0)
\, P_{\pm,k}(\lambda,k_0)^*, \quad k\in\bbN_0, \lb{4.25}
\end{align} 
where $\{P_{\pm,k}(\lambda,k_0)\}_{k\in\bbN_0}$ is an orthonormal 
system of matrix-valued polynomials with respect to the spectral 
measure $d\nu_\pm(\lambda,k_0)$, with $P_{\pm,k}(z,k_0)$ 
of degree $k$ in $z$, $P_{\pm, 0}(z,k_0)=I_m$. One verifies, 
\begin{align}
P_{+,k}(z,k_0)&=\phi(z,k_0+k,k_0-1), \lb{4.25a} \\
P_{-,k}(z,k_0)&=\theta(z,k_0-k,k_0), \quad k\in\{-1\}\cup\bbN_0,
\,\, k_0\in\bbZ, \,\, z\in\bbC, \lb{4.25b}
\end{align}
with $\phi(z,k,k_0)$ and $\theta(z,k,k_0)$ defined in 
\eqref{4.2}, \eqref{4.3}. 
\medskip

Given these preliminaries and introducing the diagonal Green's
matrix by 
\begin{equation}\lb{4.35}
g(z,k)=G(z,k,k), \quad z\in\bbC\backslash\bbR, \,\, k\in \bbZ,
\end{equation}
we can also formulate the following uniqueness result for Jacobi
operators obtained in \cite{GKM02}.

\begin{theorem} [\cite{GKM02}] \label{t4.4} 
Assume Hypothesis\ \ref{h4.1} and let $k_0\in\bbZ$. Then any
of the following three sets of data\\ 
$(i)$ $g(z,k_0)$ and $G(z,k_0,k_0+1)$ for all $z\in\bbC_+$; \\ 
$(ii)$ $g(z,k_0)$ and $[G(z,k_0,k_0+1)+G(z,k_0+1,k_0)]$ 
for all $z\in\bbC_+$; \\
$(iii)$ $g(z,k_0)$, $g(z,k_0+1)$ for all $z\in\bbC_+$ 
and $A(k_0)$; \\
uniquely determines the matrix-valued Jacobi operator $H$ 
and hence $A=\{A(k)\}_{k\in\bbZ}$ and $B=\{B(k)\}_{k\in\bbZ}$.\\ 
\end{theorem}

The special scalar case $m=1$ of Theorem\ \ref{t4.4} is known 
and has been derived in \cite{Te98} (see also 
\cite[Sect.~2.7]{Te00}. Condition (iii)  in Theorem\ \ref{t4.4} is
specific to the Jacobi case. In the corresponding Schr\"odinger 
case, the corresponding set of data does not even uniquely
determine the potential in the scalar case $m=1$ (see, 
e.g., \cite{GS96a}, \cite{Te98}).

\section{Asymptotic Expansions of Weyl--Titchmarsh and Green's 
Matrices} \label{s3}

In this section we use Riccati-type equations to systematically
determine norm convergent expansions of Weyl--Titchmarsh and Green's 
matrices as the spectral parameter tends to infinity.

We start again with the case of Jacobi operators $H$, assuming 
Hypothesis\ \ref{h4.1}.

Insertion of the norm convergent ansatz 
\begin{equation}
m_\pm(z,k)\underset{|z|\to\infty}{=} \sum_{j=1}^\infty
m_{\pm,j}(k)z^{-j},
\quad m_{\pm,1}(z,k)=-I_m \lb{3.1}
\end{equation}
into \eqref{4.12a} and \eqref{4.12b} then yields the following recursion
relation for the coefficients $m_{\pm,j}(k)$:
\begin{align}
\begin{split}
&m_{+,1}=-I_m, \quad m_{+,2}=-B, \lb{3.2} \\
&m_{+,j+1}=B m_{+,j}
-\sum_{\ell=1}^{j-1} Am_{+,j-\ell}^+Am_{+,\ell}, \quad 
j\geq2.  
\end{split}
\intertext{and}
\begin{split}
&m_{-,1}=-I_m, \quad m_{-,2}=-B, \lb{3.3} \\
&m_{-,j+1}=Bm_{-,j}
-\sum_{\ell=1}^{j-1} A^-m_{-,j-\ell}^-A^-m_{-,\ell}, \quad 
j\geq2. 
\end{split}
\end{align}
Next, rewriting \eqref{4.5f} in the form
\begin{align}
& A(k)^{-1}M_\pm(z,k+1)A(k)^{-1}M_\pm(z,k) \lb{3.4} \\
& \quad + A(k)^{-1}(zI_m-B(k+1))A(k)^{-1}M_\pm(z,k)+I_m=0, 
\quad z\in\bbC\backslash\bbR, \; k\in\bbZ, \no 
\end{align}
and inserting the norm convergent ansatz 
\begin{align}
M_+(z,k)&\underset{|z|\to\infty}{=}\sum_{j=1}^\infty M_{+,j}(k)z^{-j} 
\lb{3.5} \\
\intertext{and the asymptotic ansatz}
M_-(z,k)&\underset{|z|\to\infty}{=}-I_m z +\sum_{j=0}^\infty
M_{-,j}(k)z^{-j} \lb{3.6}
\end{align}
into \eqref{3.4} then yields the following recursion relation for the
coefficients $M_{\pm,j}(k)$:
\begin{align}
\begin{split}
& M_{+,1}=-A^2, \quad M_{+,2}=-AB^+A, \lb{3.7} \\
& M_{+,j+1}=AB^+A^{-1}M_{+,j}
-\sum_{\ell=1}^{j-1} AM^+_{+,j-\ell}A^{-1}M_{+,\ell}, \quad j\geq 2
\end{split}
\intertext{and}
\begin{split}
& M_{-,0}=B, \quad M_{-,1}=(A^-)^2, \quad M_{-,2}=A^-B^-A^-, \lb{3.8} \\
& M_{-,j+1}=-B(A^-)^{-1}M_{-,j}^-A^-
+\sum_{\ell=0}^{j} M_{-,j-\ell}(A^-)^{-1}M_{-,\ell}^-A^-, \quad j\geq 2.
\end{split}
\end{align}

\begin{remark} \label{r3.1} 
In the continuous cases of Schr\"odinger and Dirac-type operators
discussed in detail in \cite{CG01} and \cite{CG02}, establishing the
existence of appropriate asymptotic expansions was a highly nontrivial
endeavor. Here in the discrete context, hypothesis \eqref{4.1a} together
with \eqref{4.10} immediately yields a norm convergent expansion of
$m_\pm(z,k)$ as $|z|\to\infty$. By \eqref{4.11} and \eqref{4.12}, this
immediately yields the existence of an asymptotic expansion for
$M_\pm(z,k)$ as $|z|\to\infty$. By inspection, these expansions are of
the form \eqref{3.1}, \eqref{3.5}, and \eqref{3.6}.
\end{remark}

Given the asymptotic expansions \eqref{3.5} and \eqref{3.6} for
$M_{\pm}(z,k)$ as $|z|\to\infty$, one can of course derive analogous
asymptotic expansions for the $2m\times 2m$ Weyl--Titchmarsh matrix
$M(z,k)$ in \eqref{4.21}--\eqref{4.23b}. For the $(1,1)$-block matrix
element of $M(z,k)$ one obtains the norm convergent expansion for $|z|$
sufficiently large,
\begin{align}
&g(z,k)=G(z,k,k)=M_{1,1}(z,k)
\underset{|z|\to\infty}{=}\sum_{j=1}^\infty r_j(k)z^{-j}
\lb{3.9} \\
& \quad \, \underset{|z|\to\infty}{=} -I_m z^{-1} -B(k)z^{-2}
-[A(k-1)^2+A(k)^2+B(k)^2]z^{-3} \no \\
& \qquad\qquad -[B(k)^3+A(k-1)B(k-1)A(k-1)+A(k)B(k+1)A(k) \no
\\ &\qquad\qquad +B(k)A(k)^2+B(k)A(k-1)^2+A(k)^2B(k)+A(k-1)^2B(k)]z^{-4}
\no \\ 
& \qquad\qquad +\Oh(z^{-5}), \quad k\in\bbZ.  \lb{3.10}
\end{align}
Similarly,
\begin{align}
G(z,k,k+1)&\underset{|z|\to\infty}{=}-A(k)z^{-2}+O(|z|^{-3}), \quad 
k\in \bbZ, \lb{3.11} \\
G(z,k+1,k)&\underset{|z|\to\infty}{=}-A(k)z^{-2}+O(|z|^{-3}), \quad 
k\in \bbZ. \lb{3.12}
\end{align}

Next, we also recall that the $(1,1)$ and $(2,2)$-block matrices
of $M(z,k_0)$ are $m\times m$ Herglotz matrices. In particular, in
addition to \eqref{4.23c}, \eqref{4.23d}, \eqref{4.23g}, and
\eqref{4.23h} one obtains
\begin{align}
g(z,k)&=G(z,k,k)=M_{1,1}(z,k) \no \\
&=K_{1,1}(k)+\int_\bbR d\Omega_{1,1}(\lambda,k)\,
\bigg(\f1{\lambda-z}-\f{\lambda}{1+\lambda^2}\bigg) \lb{3.13} \\
&=\exp\bigg[\int_\bbR d\lambda\, \Xi (\lambda,k)
\bigg(\f1{\lambda-z}-\f{\lambda}{1+\lambda^2} \bigg) \bigg], 
\quad z\in\bbC\backslash\bbR, \; k\in\bbZ, \lb{3.14}
\end{align}
where
\begin{equation}
0\le \Xi (\dott,k)\le I_{m} \text{ a.e.} \lb{3.15}
\end{equation}
and
\begin{equation}
\Xi(\lambda,k)=\lim_{\varepsilon\downarrow
0}\f1\pi\Im(\ln(g(\lambda+i\varepsilon,k)))
\text{ for a.e.\ $\lambda\in\bbR$}.\lb{3.16}
\end{equation}
We note that the constant term in the exponent of \eqref{3.14} vanishes
because of the asymptotics \eqref{3.10}. 

\section{Borg-Type Theorems for Matrix-Valued Jacobi Operators}
\label{s4}

In this section we prove trace formulas and the discrete analog of
Borg's uniqueness theorem for matrix-valued Jacobi operators.

In the following, $\sigma(T)$ and $\sigma_{\ess}(T)$ denote the spectrum
and essential spectrum of a densely defined closed operator $T$ in a
complex separable Hilbert space.

We start with the case of matrix-valued self-adjoint Jacobi operators $H$
assuming the following hypothesis:
\begin{hypothesis} \lb{h2.1}
In addition to Hypothesis\ \ref{h4.1} suppose that $\sigma(H)\subseteq
[E_-,E_+]$ for some $-\infty < E_-<E_+<\infty$. 
\end{hypothesis}

Assuming Hypothesis\ \ref{h2.1}, trace formulas associated with Jacobi
operators then can be derived as follows. First we note that
\eqref{3.9} implies the expansion (convergent for $|z|$ sufficiently
large)
\begin{align}
&-\f{d}{dz}\ln(g(z,k))\underset{|z|\to\infty}{=}\sum_{j=1}^\infty s_j(k)
z^{-j}, \quad k\in\bbZ, \lb{2.1} \\
& \quad s_1(k)=I_m, \lb{2.2} \\
& \quad s_2(k)=B(k), \lb{2.3} \\
& \quad s_3(k)=2A(k-1)^2+2A(k)^2+B(k)^2, \;\, \text{etc.} \lb{2.4}
\end{align}

Moreover, by \eqref{3.14},
\begin{align}
\f{d}{dz}\ln(g(z,k))&=\int_{E_-}^{E^+} d\lambda\, (\lambda-z)^{-2} \, 
\Xi(\lambda,k) + \int_{E_+}^\infty d\lambda \, (\lambda-z)^{-2} I_m \no
\\
&=(E_+-z)^{-1} I_m + \int_{E_-}^{E^+} d\lambda\, (\lambda-z)^{-2} \, 
\Xi(\lambda,k),  \lb{2.5}
\end{align}
where we used
\begin{equation}
\Xi(\lambda,k)=\begin{cases} 0, & \lambda<E_-, \\
I_m, & \lambda>E_+. \end{cases}  \lb{2.6}
\end{equation}

\begin{theorem}  \lb{thm2.2}
Assume Hypothesis\ \ref{h2.1}. Then $($cf.\ \eqref{2.1}$)$, 
\begin{equation}
s_j(k)=\f{1}{2}(E_-^{j-1}+E_+^{j-1})I_m+\f{1}{2}(j-1)\int_{E_-}^{E_+}
d\lambda\, \lambda^{j-2}[I_m-2\,\Xi(\lambda,k)], \quad j\in\bbN,
\; k\in\bbZ. \lb{2.7}
\end{equation}
Explicitly, for all $k\in\bbZ$,
\begin{align}
& B(k)=\f{1}{2}(E_-+E_+) I_m + \f{1}{2}\int_{E_-}^{E_+} d\lambda \, 
[I_m-2\,\Xi(\lambda,k)], \lb{2.8} \\
& 2A(k-1)^2+2A(k)^2+B(k)^2=\f{1}{2}(E_-^2 +E_+^2) I_m +
\int_{E_-}^{E_+} d\lambda \,  \lambda[I_m-2\,\Xi(\lambda,k)], \;\,
\text{etc.} \lb{2.9}
\end{align}
\end{theorem}
\begin{proof}
By \eqref{2.1} and \eqref{2.5} one infers
\begin{align}
-\f{d}{dz}\ln(g(z,k&))=\f{1}{2}\bigg[\f{1}{z-E_+}+\f{1}{z-E_-}\bigg]I_m
+\f{1}{2}\int_{E_-}^{E_+} d\lambda\, (\lambda-z)^{-2}
[I_m-2\,\Xi(\lambda,k)]  \lb{2.10} \\
& \underset{|z|\to \infty}{=} \sum_{j=1}^\infty s_j(k) z^{-j}. \lb{2.11}
\end{align}
Expanding \eqref{2.10} with respect to $z^{-1}$ and comparing powers of
$z^{-j}$ then yields \eqref{2.7}. Equations \eqref{2.8} and \eqref{2.9}
are then clear from \eqref{2.3} and \eqref{2.4}.
\end{proof}

\begin{remark} \lb{r2.3}
In the scalar case $m=1$, Theorem\ \ref{thm2.2} was 
first derived in \cite{GS96} assuming $A(k)=1$, $k\in\bbZ$. The 
Jacobi case for half-lines was explicitly discussed in this vein 
in \cite{GS97}. The trace formula \eqref{2.8} for full-line Jacobi
operators in the case $m=1$ (and other trace formulas) can be found in
\cite{Te98}, \cite[Sect.\ 6.2]{Te00}. The current matrix-valued trace
formula for $m\geq 2$ is new.
\end{remark}

Next we turn to a Borg-type theorem for matrix-valued Jacobi operators.
To set the stage we first recall Flaschka's result \cite{Fl75} for scalar
(i.e., $m=1$) periodic Jacobi operators, the direct analog of Borg's
theorem for periodic one-dimensional Schr\"odinger operators
originally proved in \cite{Bo46}. 

\begin{theorem}  $($\cite{Fl75}$)$ \lb{t2.4}
Suppose $a$ and $b$ are periodic real-valued sequences in
$\ell^\infty(\bbZ)$ with the same period and $a(k)>0$, $k\in\bbZ$. Let
$h(a,b)=aS^+ +a^-S^- +b$ be the associated self-adjoint  Jacobi operator
in $\ell^2(\bbZ)$ $($cf.~\eqref{4.8b} for $m=1)$ and suppose that
$\sigma(h(a,b))=[E_-,E_+]$ for some $E_-<E_+$.  Then,
\begin{equation}
a(k)=(E_+-E_-)/4 , \quad  b(k)=(E_-+E_+)/2, \quad k\in\bbZ.  \lb{2.12}
\end{equation}
\end{theorem}

While uniqueness results such as Theorem~\ref{t2.4} are described as
Borg-type theorems, other types of inverse spectral results are also
described as Borg-type theorems as mentioned for Schr\"odinger
operators in Remark~\ref{r1.1}. We also note that
Theorem\ \ref{t2.4} is quite different from a recent result by Killip
and Simon \cite[Thm.\ 10.1]{KS03} which states that $a(k)=1$, $k\in\bbZ$
and $\sigma(h)\subseteq [-2,2]$ implies $b(k)=0$, $k\in\bbZ$.

As shown in \cite{CG02} and \cite{CGHL00}, periodicity is not the key
ingredient in Borg-type theorems such as Theorem\ \ref{t2.4}. In
fact, it was shown there that the more general notion of being
reflectionless is sufficient for Borg-type theorems to hold and we will
turn to this circle of ideas next. We note that the class of
reflectionless interactions include periodic and certain cases of
quasi-periodic and almost-periodic interactions.

Following \cite{CG02} and \cite{CGHL00} , we now define reflectionless
matrix-valued Jacobi operators as follows:  

\begin{definition}\lb{d2.5}
Assume Hypothesis\ \ref{h2.1}. Then the matrix-valued coefficients 
$A, B$ are called {\it reflectionless} if for all $k\in\bbZ$,
\begin{equation}
\Xi(\lambda,k)=\f12 I_m 
\text{  for a.e.\ $\lambda\in\sigma_{\ess}(H)$} \lb{2.13}
\end{equation}
with $\Xi(\cdot,k)$ defined in \eqref{3.16}.
\end{definition}

Since hardly any confusion can arise, we will also call $H=H(A,B)$ 
reflectionless if \eqref{2.13} is satisfied. 

\begin{remark}\lb{r2.6}
In the next theorem we will prove an inverse spectral result
for matrix-valued Jacobi operators $H(A,B)$. However, in the general 
matrix-valued context, where $m\geq 2$, one cannot expect that
the spectrum of $H(A,B)$ will determine $A$ and $B$ uniquely. Indeed,
assume that $B$ is a multiple of the identity, $B(k)=b(k) I_m$ for some
$b\in\ell^\infty(\bbZ)$, $b(k)\in\bbR$, $k\in\bbZ$. In addition, let $U$
be a unitary $m\times m$ matrix and consider $\wti A(k)=UA(k)U^{-1}$,
$k\in\bbZ$. Then clearly $H(A,B)$ and $H(\wti A,B)$ are unitarily
equivalent and hence the spectrum of $H$ cannot uniquely determine its
coefficients. The following result, however, will illustrate a
special case where the spectrum of $H(A,B)$ does determine $A$ and $B$
uniquely. 
\end{remark}

Given Definition\ \ref{d2.5}, we now turn to a Borg-type uniqueness
theorem for $H$ and formulate the analog of Theorem\ \ref{t2.4} for
reflectionless matrix-valued Jacobi operators.

\begin{theorem}\lb{t4.7}
Assume Hypotheses\ \ref{h2.1} and suppose that $A$ and $B$ are
reflectionless. Let $H(A,B)=AS^+ +A^-S^-+B$ be the associated
self-adjoint  Jacobi operator in $\ell^2(\bbZ)^m$ $($cf.~\eqref{4.8b}$)$
and suppose that $\sigma(H(A,B))=[E_-,E_+]$ for some $E_-<E_+$. Then
$A$ and $B$ are of the form,
\begin{equation}
A(k)=\f{1}{4}(E_+-E_-) I_m, \quad  B(k)=\f{1}{2}(E_-+E_+) I_m,
\quad k\in\bbZ. \lb{2.14}
\end{equation}
\end{theorem}
\begin{proof}
By hypothesis, $\Xi(\lambda,k)=(1/2)I_m$ for a.e.\
$\lambda\in[E_-,E_+]$ and all $k\in\bbZ$. Thus the trace formula
\eqref{2.8} immediately yields \eqref{2.14} for $B$. Inserting the
formula \eqref{2.14} for $B$ into the second trace formula \eqref{2.9}
one infers
\begin{equation}
A(k-1)^2+A(k)^2=\f{1}{8}(E_+-E_-)^2 I_m, \quad k\in\bbZ. \lb{2.15}
\end{equation} 
The first order difference equation \eqref{2.15} has the solution
\begin{equation}
A(2\ell)^2=A(0)^2, \quad A(2\ell+1)^2=\f{1}{8}(E_+-E_-)^2 I_m-A(0)^2,
\quad \ell\in\bbZ. \lb{2.16}
\end{equation}
Since by hypothesis $A(0)>0$ is a self-adjoint $m\times m$ matrix, there
exists a unitary $m\times m$ matrix $U$ that diagonalizes $A(0)$ and by
\eqref{2.16}, $U$ simultaneously diagonalizes $A(k)$ for all
$k\in\bbZ$,
\begin{equation}
\wti A(k)=UA(k)U^{-1}, \quad k\in\bbZ, \lb{2.17}
\end{equation}  
where $\wti A(k)$ are diagonal matrices for all $k\in\bbZ$. By
\eqref{3.5}--\eqref{3.8} followed by an analytic continuation to all of
$\bbC_+$, $U$ also diagonalizes $M_\pm(z,k)$, 
\begin{equation}
\wti M_\pm(z,k)=UM_\pm(z,k)U^{-1}, \quad k\in\bbZ, \lb{2.18}
\end{equation}
where $\wti M_\pm(z,k)$ are diagonal matrices for all $k\in\bbZ$. The
same result of course follows from \eqref{4.5e} taking into account that
$U$ also diagonalizes $\psi_\pm(z,k,k_0)$, $\theta(z,k,k_0)$, and
$\phi(z,k,k_0)$. The resulting diagonal matrices will of course be
denoted by $\wti\psi_\pm(z,k,k_0)$, $\wti\theta(z,k,k_0)$, and
$\wti\phi(z,k,k_0)$ below.

Next, we will invoke some Herglotz function ideas. Since for all
$k\in\bbZ$,
\begin{equation}
\Xi(\lambda,k)=\begin{cases} 0, & \lambda< E_-, \\
\f{1}{2} I_m, & \lambda\in (E_-,E_+), \\ 
I_m, & \lambda> E_+, \end{cases}  \lb{2.19}
\end{equation}
we can compute $g(z,k)$ from \eqref{3.14} and obtain
\begin{align}
g(z,k)&=\exp\bigg[\int_{E_-}^\infty d\lambda\, \Xi(\lambda,k) 
\bigg(\f{1}{\lambda-z} - \f{\lambda}{1+\lambda^2}\bigg)\bigg] \no \\
&=\exp\bigg[(1/2)\int_{E_-}^{E_+} d\lambda \, 
\bigg(\f{1}{\lambda-z} - \f{\lambda}{1+\lambda^2}\bigg)I_m \no \\
& \qquad\quad\;\; + \int_{E_+}^\infty d\lambda \, 
\bigg(\f{1}{\lambda-z} - \f{\lambda}{1+\lambda^2}\bigg)I_m\bigg] \no \\
&=-C [(z-E_-)(z-E_+)]^{-1/2}I_m, \quad z\in\bbC_+, \; k\in\bbZ 
\lb{2.20}
\end{align}
for some $C>0$. However, the known asymptotic behavior \eqref{3.10} of
$g(z,k)$ as $|z|\to\infty$ then proves $C=1$ and hence
\begin{equation}
g(z,k)=-[(z-E_-)(z-E_+)]^{-1/2} I_m :=g(z), \quad z\in\bbC_+, \;
k\in\bbZ.  \lb{2.21}
\end{equation}

Next, we consider $-g(z)^{-1}$ and its Herglotz representation,
\begin{align}
-g(z)^{-1}&=M_+(z,k)-M_-(z,k)=[(z-E_-)(z-E_+)]^{1/2}I_m \no \\
&=zI_m-B+\int_{\bbR} d\Gamma(\lambda) \, (\lambda-z)^{-1} \no \\
&= \int_{\bbR} d\Omega_+(\lambda,k) \, (\lambda-z)^{-1}
+ zI_m-B+\int_{\bbR} d\Omega_-(\lambda,k) \, (\lambda-z)^{-1},  \lb{2.22}
\end{align}
where $d\Omega_\pm(\lambda,k)$ denote the measures in the Herglotz
representations of $\pm M_\pm(z,k)$,
\begin{align}
M_+(z,k)&= \int_{\bbR} d\Omega_+(\lambda,k) \, (\lambda-z)^{-1},
\lb{2.23} \\ 
-M_-(z,k)&=zI_m-B+\int_{\bbR} d\Omega_-(\lambda,k) \, (\lambda-z)^{-1}.  
\lb{2.24}
\end{align}
Actually, applying $U$ from the left and $U^{-1}$ from the right on
either side in \eqref{2.22}--\eqref{2.24}, we may replace $M_\pm$ and
$d\Gamma$, $d\Omega_\pm$ by $\wti M_\pm$ and $d\wti \Gamma$, $d\wti
\Omega_\pm$ (where in obvious notation $d\wti \Gamma$ and $d\wti
\Omega_\pm$ denote the diagonal measures associated with the diagonal
Herglotz matrices $\wti M_+-\wti M_-$ and $\wti M_\pm$, respectively).
In the following we assume that all quantities in
\eqref{2.22}--\eqref{2.24} have been replaced by their diagonal matrix
counterparts.

Next, we follow a strategy employed in \cite{GKT96} in the scalar
Jacobi case. First we note that $\Xi(\lambda,k)=(1/2)I_m$ for a.e.\
$\lambda\in (E_-,E_+)$ is equivalent to 
\begin{equation}
g(\lambda+i0)=-g(\lambda+i0)^*, \quad \lambda\in (E_-,E_+) \lb{2.25}
\end{equation}
and hence to
\begin{equation}
-g(\lambda+i0)^{-1}=[g(\lambda+i0)^{-1}]^*, \quad \lambda\in
(E_-,E_+). \lb{2.26}
\end{equation}
Here and in the following, $f(\lambda+i0)$ denotes the normal limit
$\lim_{\varepsilon\downarrow 0} f(\lambda+i\varepsilon)$. The last result
is easily seen to be equivalent to the following fact: For all
$k\in\bbZ$,
\begin{equation}
\Re(M_+(\lambda+i0,k))=\Re(M_-(\lambda+i0,k)) \, 
\text{ for a.e.\ $\lambda\in (E_-,E_+)$}. \lb{2.27}
\end{equation}
By \eqref{4.23b} and \eqref{4.23j} one computes 
\begin{align}
&UA(k)g(z,k+1)A(k)U^{-1}=\wti A(k)g(z,k+1)\wti A(k) \no \\
&\quad =UM_{\pm}(z,k)[M_-(z,k)-M_+(z,k)]^{-1}M_{\mp}(z,k)U^{-1} \no \\
&\quad =\wti M_{\pm}(z,k)[\wti M_-(z,k)-\wti M_+(z,k)]^{-1}
\wti M_{\mp}(z,k) \no \\
&\quad =\wti M_{\pm}(z,k)g(z,k)\wti M_{\mp}(z,k), \quad z\in\bbC_+, \;
k\in\bbZ. \lb{2.28}
\end{align}
Since $g(z,k)=g(z)$ is independent of $k\in\bbZ$, \eqref{2.28} implies 
\begin{equation}
\wti A(k)^2=\wti M_+(z,k)\wti M_-(z,k)=\wti M_-(z,k)\wti M_+(z,k), \quad
z\in\bbC_+, \; k\in\bbZ. \lb{2.29}
\end{equation}
Inserting $\wti M_\pm(z,k)=\Re(\wti M_\pm(z,k))+i\Im(\wti M_\pm(z,k))$
into \eqref{2.29} then explicitly yields
\begin{equation}
\Im(\wti M_+(z,k))\Re(\wti M_-(z,k))+\Re(\wti M_+(z,k))
\Im(\wti M_-(z,k))=0 \lb{2.30} 
\end{equation}
and since $\wti M_\pm(z,k)$, $\Re(\wti M_\pm(z,k))$, and 
$\Im(\wti M_\pm(z,k))$ are all diagonal matrices, we note that all
entries in \eqref{2.30} commute. Similarly, the fact
\begin{align}
& g(\lambda+i0,k_0+1) \no \\
& \quad =\wti A(k_0)^{-1}\wti M_\pm(\lambda+i0,k_0) \no \\
& \qquad \times [\wti M_-(\lambda+i0,k_0)
-\wti M_+(\lambda+i0,k_0)]^{-1}
\wti M_\mp(\lambda+i0,k_0)\wti A(k_0)^{-1} \no \\ 
& \quad =-g(\lambda+i0,k_0+1)^*=-\wti A(k_0)^{-1}
\wti M_\mp(\lambda+i0,k_0)^* \no \\
& \qquad 
\times [\wti M_-(\lambda+i0,k_0)^*-\wti M_+(\lambda+i0,k_0)^*]^{-1}
\wti M_\pm(\lambda+i0,k_0)^*\wti A(k_0)^{-1} \lb{2.31} \\
& \hspace*{7.3cm} \text{ for a.e.\ $\lambda\in(E_-,E_+)$}  \no
\end{align}
implies
\begin{align}
&\wti M_\pm(\lambda+i0,k_0)
[\wti M_-(\lambda+i0,k_0)-\wti M_+(\lambda+i0,k_0)]^{-1}
\wti M_\mp(\lambda+i0,k_0) \no \\
& \quad =\wti M_\mp(\lambda+i0,k_0)^*
[\wti M_-(\lambda+i0,k_0)-\wti M_+(\lambda+i0,k_0)]^{-1}
\wti M_\pm(\lambda+i0,k_0)^* \lb{2.32} \\
& \hspace*{7.9cm} \text{for a.e.\ $\lambda\in(E_-,E_+)$} \no
\end{align}
since by \eqref{2.21}, 
\begin{equation}
g(z,k)=[M_-(z,k_0)-M_+(z,k_0)]^{-1}=[\wti M_-(z,k_0)-\wti
M_+(z,k_0)]^{-1} \lb{2.32a}
\end{equation}
and hence \eqref{2.25} also applies to $[\wti M_-(\lambda+i0,k_0)-\wti
M_+(\lambda+i0,k_0)]^{-1}$, $\lambda\in(E_-,E_+)$.

Inserting expression \eqref{4.4} for $\psi_\pm$ in terms of $\theta$,
$\phi$, and $M_\pm$ into \eqref{4.15} taking $\ell=k$, and inserting the
result into \eqref{2.26} then yields,
\begin{align}
&-\phi(\lambda,k,k_0)A(k_0)^{-1}M_\pm(\lambda+i0,k_0)
[M_-(\lambda+i0,k_0)-M_+(\lambda+i0,k_0)]^{-1} \no \\
& \qquad\, \times \theta(\lambda,k,k_0)^*
\no \\
&-\theta(\lambda,k,k_0)[M_-(\lambda+i0,k_0)-M_+(\lambda+i0,k_0)]^{-1}
M_\mp(\lambda+i0,k_0) A(k_0)^{-1} \no \\
& \qquad\, \times\phi(\lambda,k,k_0)^* \no \\
& \quad = -\phi(\lambda,k,k_0)A(k_0)^{-1}M_\mp(\lambda+i0,k_0)^*
[M_-(\lambda+i0,k_0)-M_+(\lambda+i0,k_0)]^{-1} \no \\
& \qquad\, \times \theta(\lambda,k,k_0)^*
\no \\
&
\qquad -\theta(\lambda,k,k_0)
[M_-(\lambda+i0,k_0)-M_+(\lambda+i0,k_0)]^{-1}
M_\pm(\lambda+i0,k_0)^* A(k_0)^{-1} \no \\
& \qquad\, \times \phi(\lambda,k,k_0)^*, \lb{2.33}
\end{align}
where we also used \eqref{2.25} and \eqref{2.32}. Applying $U$ and
$U^{-1}$ from the left and right on either side in \eqref{2.33}, using
the fact that 
\begin{equation}
\wti \phi(\ol z,k,k_0)^*=\wti \phi(z,k,k_0), \quad 
\wti \theta(\ol z,k,k_0)^*=\wti \theta(z,k,k_0), \lb{2.34}
\end{equation}
and that all diagonal matrices commute, one can rewrite \eqref{2.33} in
the form,
\begin{align}
& 2i\wti\phi(\lambda,k,k_0)\wti\theta(\lambda,k,k_0)A(k_0)^{-1}
g(\lambda+i0)  \lb{2.35} \\
& \quad \times [\Im(\wti M_-(\lambda+i0,k_0))
+ \Im(\wti M_+(\lambda+i0,k_0))]=0 \, 
\text{ for a.e.\ $\lambda\in (E_-,E_+)$}. \no 
\end{align}
Since $k\in\bbZ$ can be chosen arbitarily in \eqref{2.35}, this implies
\begin{equation}
\Im(\wti M_-(\lambda+i0,k_0))=-\Im(\wti M_+(\lambda+i0,k_0)) \, 
\text{ for a.e.\ $\lambda\in (E_-,E_+)$}. \lb{2.36}
\end{equation}
Finally, since $k_0\in\bbZ$ in \eqref{2.36} is arbitrary, one obtains
for all $k\in\bbZ$,
\begin{equation}
\Im(\wti M_-(\lambda+i0,k))=-\Im(\wti M_+(\lambda+i0,k)) \, 
\text{ for a.e.\ $\lambda\in (E_-,E_+)$} \lb{2.37}
\end{equation}
and hence also for all $k\in\bbZ$,
\begin{equation}
\Im(M_-(\lambda+i0,k))=-\Im(M_+(\lambda+i0,k)) \, 
\text{ for a.e.\ $\lambda\in (E_-,E_+)$} \lb{2.38}
\end{equation}
applying $U^{-1}$ and $U$ from the left and right on both
sides in \eqref{2.37}.  Together with \eqref{2.27} this yields for all
$k\in\bbZ$,
\begin{equation}
M_-(\lambda \mp i0,k)=M_-(\lambda \pm i0,k)^*=M_+(\lambda \pm i0,k) \, 
\text{ for a.e.\ $\lambda\in (E_-,E_+)$}. \lb{2.39}
\end{equation}
Thus, $M_-(\cdot,k)$ is the analytic continuation of $M_+(\cdot,k)$
(and vice versa) through the interval $(E_-,E_+)$. Since $d\Gamma$
is purely absolutely continuous (cf.\ \eqref{2.22}), 
\begin{equation}
d\Gamma=d\Gamma_{\rm ac}, \quad d\Gamma_{\rm pp}=d\Gamma_{\rm sc}=0
\lb{2.40} 
\end{equation}
(where $d\mu_{\rm ac}$, $d\mu_{\rm pp}$, and $d\mu_{\rm sc}$
denote the absolutely continuous, pure point, and singularly continuous
parts of a measure $d\mu$), one also infers
\begin{equation}
d\Omega_{\pm,\rm pp}=d\Omega_{\pm,\rm sc}=0. \lb{2.41}
\end{equation}
(This also follows from the fact that $M_\pm(\cdot,k)$ have analytic
continuations through $(E_-,E_+)$, see \cite[Lemma\ 5.6]{GT00}.)
Especially, \eqref{2.37} then implies the $k$-independence of 
$d\Omega_\pm(\cdot,k)$,
\begin{equation}
d\Omega_+(\lambda,k)=d\Omega_-(\lambda,k)=\f{1}{2}d\Gamma(\lambda).
\lb{2.42}
\end{equation}
Equations \eqref{3.7}, \eqref{3.8}, \eqref{2.23}, \eqref{2.24}, and
\eqref{2.42} then prove
\begin{equation}
\wti A(k)^2=\int_{E_-}^{E_+}
d\Omega_+(\lambda,k)=\f{1}{2}\int_{E_-}^{E_+}
d\Gamma(\lambda)=\int_{E_-}^{E_+} d\Omega_-(\lambda,k)=\wti A(k-1)^2,
\quad k\in\bbZ \lb{2.43}
\end{equation}
and hence also 
\begin{equation}
A(k)^2=A(k-1)^2, \quad k\in\bbZ, \lb{2.44}
\end{equation}
which proves that $A$ is independent of $k\in\bbZ$. The trace formula
\eqref{2.15} for $\wti A$,
\begin{equation}
\wti A(k-1)^2+\wti A(k)^2=\f{1}{8}(E_+-E_-)^2 I_m, \quad k\in\bbZ. 
\lb{2.45}
\end{equation}
then proves 
\begin{equation}
\wti A(k)^2=\f{1}{16}(E_+-E_-)^2 I_m, \quad k\in\bbZ. \lb{2.46} 
\end{equation}
Using \eqref{2.17} then completes the proof of \eqref{2.14}.
\end{proof}

Because of \eqref{2.42} and 
\begin{equation}
d\Gamma(\lambda)=\begin{cases} 
\f{1}{\pi}[(\lambda-E_-)(E_+-\lambda)]^{1/2}, & \lambda\in [E_-,E_+], \\
0, & \lambda \in\bbR\backslash [E_-,E_+], \lb{2.47} 
\end{cases}
\end{equation}
(cf.\ \eqref{2.22}), one obtains
\begin{equation}
\wti M_\pm(z,k)=\bigg\{-\f{1}{2}z+\f{1}{4}(E_-+E_+) \pm\f{1}{2} 
[(z-E_-)(z-E_+)]^{1/2}\bigg\}I_m, \quad z\in\bbC_+, \; k\in\bbZ.
\lb{2.48}
\end{equation}

\begin{corollary}\lb{c4.8}
Assume Hypothesis\ \ref{h2.1} in the special case $m=1$ and suppose
that
$a$ and $b$ are reflectionless. Let $h(a,b)=aS^+ +a^-S^-+b$ be the
associated self-adjoint  Jacobi operator in $\ell^2(\bbZ)$
$($cf.~\eqref{4.8b}$)$ and suppose that $\sigma(h(a,b))=[E_-,E_+]$ for
some $E_-<E_+$.  Then, 
\begin{equation}
a(k)=\f{1}{4}(E_+-E_-), \quad  b(k)=\f{1}{2}(E_-+E_+), \quad
k\in\bbZ. \lb{2.49}
\end{equation}
\end{corollary}

While Theorem\ \ref{t4.7} is new, Corollary\ \ref{c4.8} in the scalar
case $m=1$ was noted in \cite[Corollary\ 6.3]{Te98} (see also
\cite[Corollary\ 8.6]{Te00}). 

The following result can be proved in analogy to Theorems\ 4.6 and 4.8
in \cite{CGHL00}, hence we state it here without proof.

\begin{theorem} \lb{t4.9} 
In addition to Hypothesis\ \ref{h2.1}, suppose that $A$ and
$B$ are periodic with the same period. Let $H(A,B)=AS^+ +A^-S^-+B$ be the
associated self-adjoint  Jacobi operator in $\ell^2(\bbZ)^m$ and suppose
that $H(A,B)$ has uniform spectral multiplicity $2m$. Then $H(A,B)$ is
reflectionless and hence for all $k\in\bbZ$,
\begin{equation}
\Xi(\lambda,k)=\f{1}{2} I_m \, \text{ for a.e.\
$\lambda\in\sigma_{\ess}(H(A,B))$.}  \lb{2.50}
\end{equation}
In particular, assume that $A$ and $B$ are periodic with the same period,
that $H(A,B)$ has uniform spectral multiplicity $2m$, and that 
$\sigma(H(A,B))=[E_-,E_+]$ for some $E_-<E_+$.  Then
$A$ and $B$ are of the form,
\begin{equation}
A(k)=\f{1}{4}(E_+-E_-) I_m, \quad  B(k)=\f{1}{2}(E_-+E_+) I_m,
\quad k\in\bbZ. \lb{2.51}
\end{equation}
\end{theorem}

In connection with the special case $m=1$ in Theorem \ref{t2.4} we note
that scalar Jacobi operators automatically have uniform spectral
multiplicity $2$ since the product of the two Floquet multipliers
equals one.

\section{Borg-Type Theorems for Supersymmetric Dirac \\ Difference
Operators} \label{s5}

In our final section we turn to a Borg-type theorem for a class of
supersymmetric Dirac difference operators. Rather than developing the
results from first principles as in the case of Jacobi operators with
matrix-valued coefficients, we will employ the underlying supersymmetric
structure and reduce the case of Dirac difference operators to that of
Jacobi operators. 

We start with the following abstract result.

\begin{theorem} [\cite{De78}, \cite{GSS91}] \lb{thm4.10} 
Let $C$ be a densely defined closed operator in a separable
complex Hilbert space $\calH$ with domain $\calD(C)$ and introduce the
operator 
\begin{equation}
Q=\begin{pmatrix} 0 & C^*\\ C & 0 \end{pmatrix}, \quad
\calD(Q)=\calD(C)\oplus \calD(C^*) \lb{2.52}
\end{equation}
in $\calH\oplus\calH$. Then,
\begin{align}
&Q=Q^*, \lb{2.53} \\
&Q^2=\begin{pmatrix} C^*C & 0 \\ 0 & CC^* \end{pmatrix}, \lb{2.54} \\
&\sigma_3 Q\sigma_3 =-Q, \quad 
\sigma_3=\begin{pmatrix} I_\calH & 0 \\ 0 & -I_\calH \end{pmatrix}, 
\lb{2.55} \\
& (Q-zI_{\calH\oplus\calH})^{-1}=\begin{pmatrix} 
z(C^*C-z^2I_\calH)^{-1} & C^*(CC^*-z^2I_\calH)^{-1} \\
C(C^*C-z^2I_\calH)^{-1} & z(CC^*-z^2I_\calH)^{-1}\end{pmatrix},
\lb{2.56}\\
& \hspace*{5.3cm} z^2\in \bbC\backslash\{\sigma(C^*C)\cup \sigma(CC^*)\}.
\no
\end{align}
\end{theorem}

\medskip

In addition, we mention the following facts: 
\begin{align}
&I_\calH+\zeta (CC^*-\zeta I_\calH)^{-1}\supseteq 
C(C^*C-\zeta I\calH)^{-1}C^*, \lb{2.57} \\
& \hspace*{2.88cm}  \zeta\in
\bbC\backslash\{\sigma(C^*C)\cup \sigma(CC^*)\}, \no \\
&I_\calH+\zeta (C^*C-\zeta I_\calH)^{-1}\supseteq 
C^*(CC^*-\zeta I\calH)^{-1}C, \lb{2.58} \\
& \hspace*{2.88cm} \zeta\in \bbC\backslash\{\sigma(C^*C)\cup
\sigma(CC^*)\}. \no
\end{align}
Moreover,
\begin{equation}
QU(z)=zU(z), \quad U(z)=\begin{pmatrix} u_1(z) \\
u_2(z) \end{pmatrix}  \lb{2.58a}
\end{equation} 
implies 
\begin{equation}
C^*u_2(z)=zu_1(z), \quad Cu_1(z)=zu_2(z) \lb{2.58b}
\end{equation}
and hence
\begin{equation}
C^*Cu_1(z)=z^2u_1(z), \quad CC^*u_2(z)=z^2u_2(z). 
\lb{2.58c}
\end{equation}
Conversely,
\begin{equation}
C^*Cu_1(z)=z^2u_1(z), \quad z\neq 0, \lb{2.58d}
\end{equation}
implies 
\begin{equation}
QU(z)=zU(z), \quad U(z)=\begin{pmatrix} u_1(z) \\
(1/z)Cu_1(z) \end{pmatrix}  \lb{2.58e}
\end{equation} 
and 
\begin{equation}
CC^*u_2(z)=z^2u_2(z), \quad z\neq 0, \lb{2.58f}
\end{equation}
implies 
\begin{equation}
QU(z)=zU(z), \quad U(z)=\begin{pmatrix} (1/z)C^*u_2(z) \\
u_2(z) \end{pmatrix}.  \lb{2.58g}
\end{equation} 

In order to apply this setup to finite difference Dirac-type systems
(cf.\ \cite{CG03}), we introduce the following hypothesis.

\begin{hypothesis} \lb{h4.11}
Let $m\in\bbN$ and consider the sequence of invertible
$m\times m$ matrices 
\begin{align}
& \rho=\{\rho(k)\}_{k\in\bbZ}\in \bbC(\bbZ)^{m\times m}, \quad 
\rho(k)=\rho(k)^*, \; k\in\bbZ,
\lb{2.59} \\
& \chi=\{\chi(k)\}_{k\in\bbZ}\in \bbC(\bbZ)^{m\times m}, 
\lb{2.60} \\
& {\det}_{\bbC^m}(\rho(k))\neq 0, \quad {\det}_{\bbC^m}(\chi(k))\neq 0,
\quad  k\in\bbZ. \lb{2.61}
\end{align}
In addition, we assume that $\rho(k)$ is a  positive definite diagonal
$m\times m$ matrix 
\begin{equation}
\rho(k)=\diag(\rho_1(k),\dots,\rho_m(k)), \quad \rho_j(k)>0, \; 1\leq
j\leq m, \; k\in\bbZ, \lb{2.62}
\end{equation} 
and that $\rho(k)$ and $\chi(k)$ are uniformly bounded
with  respect to $k\in\bbZ$, that is, there exists a $C>0$,
such that 
\begin{equation}
\|\rho(k)\|_{\bbC^{m\times m}}+\|\chi(k)\|_{\bbC^{m\times m}}\leq C,
\quad k\in\bbZ. \lb{2.63}
\end{equation}
Finally, we suppose that $\rho\chi^+$ and $\chi\rho$ are positive
definite, 
\begin{equation}
\rho(k)\chi(k+1)>0, \;\; \chi(k)\rho(k)>0, \; k\in\bbZ. \lb{2.64}
\end{equation}
\end{hypothesis}

\medskip

Assuming Hypothesis\ \ref{h4.11}, we thus introduce the bounded linear
operator
\begin{equation}
E=\rho^- S^- + \chi, \quad \calD(E)=\ell^2(\bbZ)^m, \lb{2.65}
\end{equation}
on $\ell^2(\bbZ)^m$ and the bounded Dirac-type difference operator
\begin{equation}
D=\begin{pmatrix} 0 & E^* \\ E & 0 \end{pmatrix} =
S_\rho+X, \quad \calD(D)=\ell^2(\bbZ)^m\oplus\ell^2(\bbZ)^m  \lb{2.66}
\end{equation}
on $\ell^2(\bbZ)^m\oplus\ell^2(\bbZ)^m$, where
\begin{equation}
S_\rho=\begin{pmatrix} 0 & \rho S^+ \\ 
\rho^- S^- & 0 \end{pmatrix}, \quad X= \begin{pmatrix} 
0 & \chi^* \\ \chi& 0 \end{pmatrix}. \lb{2.67}
\end{equation}
One then computes
\begin{align}
H_1&=E^*E=A_1S^+ +A_1^-S^-+B_1, \lb{2.68} \\
H_2&=EE^*=A_2S^+ +A_2^-S^-+B_2, \lb{2.69}  
\end{align}
where
\begin{align}
A_1(k)&=\rho(k)\chi^+(k)>0, \quad
B_1(k)=(\rho(k))^2+\chi(k)^*\chi(k), \; k\in\bbZ, \lb{2.70} \\ 
A_2(k)&=\chi(k)\rho(k)>0, \quad
B_2(k)=(\rho(k)^-)^2+\chi(k)\chi(k)^* \lb{2.71} 
\end{align}
and notes that $H_1\geq 0$ and $H_2\geq 0$ are matrix-valued Jacobi
operators in $\ell^2(\bbZ)^m$ of the form \eqref{4.8b}. 

\begin{remark} \lb{r4.12a}
We note that $D$ with a positive definite diagonal $m\times m$ matrix
$\rho$ in \eqref{2.66} represents a normal form of Dirac-type
difference in the following sense: Assume Hypothesis\ \ref{h4.11} except
for the condition that $\rho(k)$ is a positive definite diagonal matrix
for all $k\in\bbZ$. Then, following \cite[Lemma\ 2.3]{CG03}, there
exists a sequence of unitary matrices
$U(\rho)=\{U(\rho,k)\}_{k\in\bbZ}\in\bbC(\bbZ)^{2m\times 2m}$
such that
\begin{equation}
U(\rho)(S_\rho+X)U(\rho)^{-1}=S_{\widehat\rho}+\widehat X,
\lb{2.71a}
\end{equation}
where $\widehat\rho$ is diagonal and positive definite and
$\widehat X$ is of the form
\begin{equation}
\widehat X=U(\rho)XU(\rho)^{-1}=
\begin{pmatrix} 0 & {\widehat\chi}^* \\ \widehat\chi & 0
\end{pmatrix}, 
\lb{2.71b}
\end{equation}
with $\widehat\chi\in\bbC(\bbZ)^{m\times m}$. Indeed, denote by 
$Q_\rho (k)\in\bbC^{m\times m}$ a unitary matrix such that
$Q_\rho (k)\rho(k)Q_\rho(k)^{-1} = \wti \rho(k)$, where 
$\wti \rho(k)\in\bbR^{m\times m}$ is diagonal and self-adjoint for all
$k\in\bbZ$. Then,
\begin{equation}
{U_{\rho}}(S_{\rho} +X)U_{\rho}^{-1} = S_{\wti \rho}+ \wti X,
\quad \wti X = {U_{\rho}}XU_{\rho}^{-1}, \quad
U_{\rho}=\begin{pmatrix}Q_\rho & 0 \\ 0 & Q_\rho^-\end{pmatrix}. 
\lb{2.71c}
\end{equation}
Next, let $\wti\varepsilon (k)\in\bbR^{m\times m}$ be
a diagonal matrix for which
$(\wti\varepsilon (k))_{\ell,\ell}\in\{+1,-1\}$,
$\ell=1,\dots,m$. Define $\varepsilon (k)\in\bbR^{m\times m}$ by
$\varepsilon (k)=\wti\varepsilon (k) \wti \varepsilon (k+1)$ and
choose $\wti\varepsilon (k)$ so that 
$\widehat\rho=\varepsilon \wti \rho >0$. Then,
\begin{equation}
U_{\varepsilon }(S_{\wti \rho} + \wti X)U_\epsilon^{-1}  =
S_{\widehat \rho} +\widehat X, \quad 
\widehat X = U_{\varepsilon}
\wti X U_\varepsilon^{-1}, \quad U_{\varepsilon}
= \begin{pmatrix}\wti\varepsilon & 0 \\ 0 & \wti\varepsilon 
\end{pmatrix}. \lb{2.71d}
\end{equation}
Thus, one obtains 
\begin{gather}
U(\rho)=U_\varepsilon U_\rho=\begin{pmatrix} 
\wti\varepsilon Q_\rho & 0 \\ 0 & \wti\varepsilon Q_\rho^-
\end{pmatrix}, \lb{2.71e} \\
\widehat X=\begin{pmatrix} 0 & {\widehat\chi}^* \\ \widehat\chi & 0
\end{pmatrix}, \quad \widehat\chi=\wti\varepsilon Q^-_\rho\chi
Q_\rho^{-1}\wti\varepsilon. \lb{2.71f}
\end{gather}
\end{remark}  

Next, we introduce the $2m\times 2m$ Weyl--Titchmarsh matrices
associated with $D$ (cf.\ \cite{CG03}) by 
\begin{equation}
M^D(z,k)=\big(M^D_{j,j^\prime}(z,k)\big)_{j,j^\prime=1,2}\, , 
\quad z\in\bbC\backslash\bbR, \,\, k\in\bbZ, \lb{2.72}
\end{equation}
where
\begin{align}
M^D_{1,1}(z,k)&=\big[M^D_{-}(z,k)-M^D_{+}(z,k)\big]^{-1}, \lb{2.73}
\\ M^D_{1,2}(z,k)&=2^{-1}\big[M^D_{-}(z,k)-M^D_{+}(z,k)\big]^{-1}
\big[M^D_{-}(z,k)+M^D_{+}(z,k)\big], \lb{2.74} \\
M^D_{2,1}(z,k)&=2^{-1}\big[M^D_{-}(z,k)+M^D_{+}(z,k)\big]
\big[M^D_{-}(z,k)-M^D_{+}(z,k)\big]^{-1}, \lb{2.75} \\
M^D_{2,2}(z,k)&=M^D_{\pm}(z,k)
\big[M^D_{-}(z,k)-M^D_{+}(z,k)\big]^{-1}M_{\mp}(z,k), \lb{2.76} 
\end{align}
and similarly, the $2m\times 2m$ Weyl--Titchmarsh matrices
associated with the Jacobi operators $H_\ell$, $\ell=1,2$ (cf.\
\eqref{4.21}--\eqref{4.23b}) by 
\begin{equation}
M^{H_\ell}(z,k)
=\big(M^{H_\ell}_{j,j^\prime}(z,k)\big)_{j,j^\prime=1,2}\,, 
\quad \ell=1,2, \; z\in\bbC\backslash\bbR, \; k\in\bbZ, \lb{2.77}
\end{equation}
where
\begin{align}
M^{H_\ell}_{1,1}(z,k)
&=\big[M^{H_\ell}_{-}(z,k)-M^{H_\ell}_{+}(z,k)\big]^{-1}, \lb{2.78}
\\
M^{H_\ell}_{1,2}(z,k)&=2^{-1}\big[M^{H_\ell}_{-}(z,k)
-M^{H_\ell}_{+}(z,k)\big]^{-1}
\big[M^{H_\ell}_{-}(z,k)+M^{H_\ell}_{+}(z,k)\big], \lb{2.79} \\
M^{H_\ell}_{2,1}(z,k)&=2^{-1}\big[M^{H_\ell}_{-}(z,k)
+M^{H_\ell}_{+}(z,k)\big]
\big[M^{H_\ell}_{-}(z,k)-M^{H_\ell}_{+}(z,k)\big]^{-1}, \lb{2.80} \\
M^{H_\ell}_{2,2}(z,k)&=M^{H_\ell}_{\pm}(z,k)
\big[M^{H_\ell}_{-}(z,k)-M^{H_\ell}_{+}(z,k)\big]^{-1}M_{\mp}(z,k).
\lb{2.81} 
\end{align} 
The supersymmetric formalism \eqref{2.52}--\eqref{2.56} then
implies the following relations between $M^D(z)$ and $M^{H_\ell}(z)$,
$\ell=1,2$.

\begin{theorem} \lb{t4.12}
Assume Hypothesis\ \ref{h4.11} and let $z\in\bbC\backslash\sigma(D)$,
$k\in\bbZ$. Then, 
\begin{align}
M^D_\pm(z,k)&=-z^{-1}\rho(k)+z^{-1}\rho(k)^{-1/2}M^{H_1}_\pm(z^2,k)
\rho(k)^{-1/2}, \lb{2.82} \\
M^D_\pm(z,k)&=-z\rho(k)^{-1}-z\rho(k)^{-1/2}
\big[\chi(k)^{-1}M^{H_2}_\pm
(\chi(k)^{-1})^*\big]^{-1}\rho(k)^{-1/2}. \lb{2.83}
\end{align}
\end{theorem}
\begin{proof}
We freely employ the relations \eqref{2.58a}--\eqref{2.58g} and some 
results from \cite{CG03}. Let 
\begin{equation}
U_\pm(z,k,k_0)=\begin{pmatrix} u_{1,\pm}(z,k,k_0) \\ 
u_{2,\pm} (z,k,k_0) \end{pmatrix}, \quad u_{1,\pm}(z,k_0,k_0)=I_m 
\lb{2.84}
\end{equation}
be the normalized Weyl--Titchmarsh solutions associated with $D$. Then
by equations (2.35a)--(2.35c) and (2.96) in \cite{CG03}, 
\begin{equation}
u_{2,\pm}(z,k+1,k_0)u_{1,\pm}(z,k,k_0)^{-1}
=-\rho(k)^{-1/2}M^D_\pm(z,k)\rho(k)^{1/2}. \lb{2.85} 
\end{equation}
Moreover, using $(C^*u_{2,\pm})(z,k_0,k_0)=zu_{1,\pm}(z,k_0,k_0)=zI_m$
(cf.\ \eqref{2.58b}), one derives 
\begin{equation}
u_{2,\pm}(z,k_0,k_0)=(\chi(k_0)^*)^{-1}
\big[zI_m+\rho(k_0)^{1/2}M^D_\pm(z,k_0)\rho(k_0)^{1/2}\big]. \lb{2.86}
\end{equation}
Next, let 
\begin{equation}
\psi_{\ell,\pm}(z,k,k_0), \quad \psi_{\ell,\pm}(z,k_0,k_0)=I_m 
\lb{2.87}
\end{equation}
be the normalized Weyl--Titchmarsh solutions associated with $H_\ell$,
$\ell=1,2$. Then by \eqref{4.5e},
\begin{equation}
M^{H_\ell}_\pm(z,k)=-A_\ell(k)\psi_{\ell,\pm}(z,k+1,k_0)
\psi_{\ell,\pm}(z,k,k_0)^{-1}. \lb{2.88}
\end{equation}
Given the uniqueness of Weyl--Titchmarsh solutions for $D$ and
$H_\ell$, $\ell=1,2$, \eqref{2.58a}--\eqref{2.58g} yield
\begin{align}
u_{1,\pm}(z,k,k_0)&=\psi_{1,\pm}(z^2,k,k_0), \lb{2.89} \\
u_{2,\pm}(z,k,k_0)&=(1/z)Eu_{1,\pm}(z,k,k_0) \no \\
&=(1/z)[\rho^-(k)u_{1,\pm}(z,k-1,k_0)+\chi(k)u_{1,\pm}(z,k,k_0)] \no \\
&=(1/z)[\rho^-(k)\psi_{1,\pm}(z^2,k-1,k_0)
+\chi(k)\psi_{1,\pm}(z^2,k,k_0)]. \lb{2.90} 
\end{align}
Thus,
\begin{align}
&u_{2,\pm}(z,k+1,k_0)u_{1,\pm}(z,k,k_0)^{-1}
=-\rho(k)^{-1/2}M^D_\pm(z,k)\rho(k)^{1/2} \no \\
& \quad =(1/z)\big[\rho(k)\psi_{1,\pm}(z^2,k,k_0)+\chi(k+1)
\psi_{1,\pm}(z^2,k+1,k_0) \big]\psi_{1,\pm}(z^2,k,k_0)^{-1} \no \\
& \quad = (1/z)\big[\rho(k)+\chi(k+1)(-A_1(k))^{-1}
M^{H_1}_\pm(z^2,k)\big] \no \\
& \quad = (1/z)\big[\rho(k)-\rho(k)^{-1}M^{H_1}_\pm(z^2,k) \big] 
\lb{2.91}
\end{align}
proving \eqref{2.82}.

Similarly, the uniqueness of Weyl--Titchmarsh solutions also yields
\begin{equation}
u_{2,\pm}(z,k,k_0)=\psi_{2,\pm}(z^2,k,k_0)d_\pm(z,k_0) \lb{2.92}
\end{equation}
for some constant $m\times m$ matrix $d_\pm(z,k_0)$. Thus,
\begin{equation}
d_\pm(z,k_0)=u_{2,\pm}(z,k_0,k_0)=(\chi(k_0)^*)^{-1}[zI_m+
\rho(k_0)^{1/2}M^D_\pm(z,k_0)\rho(k_0)^{1/2}]. \lb{2.93}
\end{equation}
One then computes,
\begin{align}
&u_{2,\pm}(z,k_0+1,k_0)=-\rho(k_0)^{-1/2}M^D_\pm(z,k_0)\rho(k_0)^{1/2}
\no \\
& \quad =\psi_{2,\pm}(z^2,k_0+1,k_0)(\chi(k_0)^*)^{-1}[zI_m+
\rho(k_0)^{1/2}M^D_\pm(z,k_0)\rho(k_0)^{1/2}] \no \\
& \quad =-(A_2(k_0))^{-1}M^{H_2}_\pm(z^2,k_0)(\chi(k_0)^*)^{-1}[zI_m+
\rho(k_0)^{1/2}M^D_\pm(z,k_0)\rho(k_0)^{1/2}] \no \\
& \quad = -\rho(k_0)^{-1}\chi(k_0)^{-1}M^{H_2}_\pm(z^2,k_0)
(\chi(k_0)^*)^{-1}[zI_m+
\rho(k_0)^{1/2}M^D_\pm(z,k_0)\rho(k_0)^{1/2}]. \lb{2.94}
\end{align}
Hence,
\begin{align}
&\rho(k_0)^{1/2}M^D_\pm(z,k_0)\rho(k_0)^{1/2} \no \\
&\quad =\chi(k_0)^{-1}M^{H_2}_\pm(z^2,k_0)
(\chi(k_0)^*)^{-1}[zI_m+
\rho(k_0)^{1/2}M^D_\pm(z,k_0)\rho(k_0)^{1/2}] \lb{2.95}
\end{align}
and since $k_0\in\bbZ$ is arbitrary,
\begin{align}
&\rho(k)^{1/2}M^D_\pm(z,k)\rho(k)^{1/2} \lb{2.96} \\
&\quad =\chi(k)^{-1}M^{H_2}_\pm(z^2,k)
(\chi(k)^*)^{-1}[zI_m+
\rho(k)^{1/2}M^D_\pm(z,k)\rho(k)^{1/2}], \quad k\in\bbZ. \no
\end{align}
Solving \eqref{2.96} for $M^D_\pm(z,k)$ then yields \eqref{2.83}.
\end{proof}

According to \eqref{4.23d}, $M^D(z,k)$ is a matrix-valued Herglotz
function of  rank $2m$ with exponential Herglotz representation
\begin{align}
M^D(z,k)&=\exp\bigg[C^D(k)+\int_\bbR d\lambda\, 
\Upsilon^D (\lambda,k)
\bigg(\f1{\lambda-z}-\f{\lambda}{1+\lambda^2} \bigg) \bigg], 
\lb{2.100}
\end{align}
where
\begin{align}
C^D(k)&=C^D(k)^*, \quad 0\le \Upsilon^D (\dott,k)\le I_{2m} 
\text{ a.e.,} \lb{2.101} \\
\Upsilon^D (\lambda,k)&=\lim_{\varepsilon\downarrow
0}\f1\pi\Im(\ln(M(\lambda+i\varepsilon,k)))
\text{ for a.e.\ $\lambda\in\bbR$}.\lb{2.102}
\end{align}

Following \cite{CG02} we now define reflectionless
matrix-valued Dirac-type operators as follows:  

\begin{definition}\lb{d4.13}
Assume Hypothesis\ \ref{h4.11}. Then the matrix-valued coefficients 
$\rho, \chi$ are called {\it reflectionless} if for all $k\in\bbZ$,
\begin{equation}
\Upsilon^D (\lambda,k)=\f12 I_{2m} 
\text{  for a.e.\ $\lambda\in\sigma_{\ess}(D)$}. \lb{2.103}
\end{equation}
\end{definition}

We also call $D=D(\rho,\chi)$  reflectionless if \eqref{2.103}
holds. 

\begin{remark} \lb{r4.14}
The definition \eqref{2.13} of reflectionless Jacobi operators $H$ 
and the definition \eqref{2.103} of reflectionless Dirac operators
$D$ can be replaced by the more stringent requirement that for all
$k\in\bbZ$,
\begin{equation}
M_+(\lambda + i0,k)=M_-(\lambda - i0,k) \, 
\text{ for a.e.\ $\lambda\in \sigma_{\ess}(H)$} \lb{2.104}
\end{equation}
for Jacobi operators $H$ and similarly for Dirac-type operators
$D$ replacing $M_\pm(z,k)$ by $M^D_\pm(z,k)$. This yields a unified
definition of the notion of reflectionless matrix-valued Jacobi and
Dirac-type operators. It is easy to see that \eqref{2.104} implies
\eqref{2.13}. The converse is more subtle and was proved by Sodin
and Yuditskii in \cite{SY95} and \cite{SY96} for scalar
Schr\"odinger and Jacobi operators $H$ under the assumption that
$\sigma(H)$ is a homogeneous set. Their proof extends to the
present matrix-valued setting. 

In the special case of a Borg-type situation with
$\sigma(H)=[E_-,E_+]$, we explicitly derived \eqref{2.104} in
\eqref{2.39}. In this particular case, $M_-(\cdot,k)$ is the
analytic continuation of $M_+(\cdot,k)$ (and vice versa) through the
interval $(E_-,E_+)$. In general, the homogeneous set $\sigma(H)$
may be a Cantor set (of positive Lebesgue measure) and then
$M_\pm(\cdot,k)$ are called {\it pseudo-continuable} through
$\sigma(H)$.
\end{remark}

\begin{lemma} \lb{l4.15}
Assume Hypothesis\ \ref{h4.11}. If $D$ is reflectionless $($in the
sense of \eqref{2.103}$)$ then $H_\ell$, $\ell=1,2$, are
reflectionless  $($in the sense of \eqref{2.13}$)$.
\end{lemma}
\begin{proof}
$D$ being reflectionless in the sense of \eqref{2.103} is equivalent
to  the assertion that for all $k\in\bbZ$, 
$M^D(\lambda+i0,k)$ is skew-adjoint for a.e.\
$\lambda\in\sigma_{\ess}(D)$. Equivalently, for all
$k\in\bbZ$,
\begin{equation}
M^D(\lambda+ i0,k)=iC(\lambda,k) \,
\text{ with $C(\lambda,k)=C(\lambda,k)^*$ for a.e.\
$\lambda\in\sigma_{\ess}(D)$}.  \lb{2.105}
\end{equation}
In fact, $C(\lambda,k)=\Im(M^D(\lambda+ i0,k))\geq 0$. As a
consequence, also all block submatrices of $M^D(\lambda+ i0,k)$,
symmetric with respect to the diagonal of $M^D(\lambda+ i0,k)$, are
skew-adjoint. In particular, the two $m\times m$ diagonal blocks
of $M^D(\lambda+i0,k)$ satisfy for all $k\in\bbZ$,
\begin{align}
M^D_{\ell,\ell}(\lambda+ i0,k)&=iC_{\ell,\ell}(\lambda,k) \,
\text{ with $C_{\ell,\ell}(\lambda,k)\geq 0$ for a.e.\
$\lambda\in\sigma_{\ess}(D)$}.  \lb{2.106}
\end{align}
Equation \eqref{2.82} implies 
\begin{equation}
[M^D_-(z,k)-M^D_+(z,k)]^{-1}  
=\rho(k)^{1/2} z
\big[M^{H_1}_-(z^2,k)-M^{H_1}_+(z^2,k)\big]^{-1}\rho(k)^{1/2}, \quad
k\in\bbZ  \lb{2.107} 
\end{equation}
and hence, for all $k\in\bbZ$ and a.e.\
$\lambda\in\sigma_{\ess}(D)$,
\begin{align}
M^D_{1,1}(\lambda+i0,k)
&=\big[M^D_-(\lambda+i0,k)-M^D_+(\lambda+i0,k)\big]^{-1} \no \\ 
&=\rho(k)^{1/2} \lambda 
\big[M^{H_1}_-(\lambda^2+i0,k)-M^{H_1}_+(\lambda^2+i0,k)\big]^{-1}
\rho(k)^{1/2}.  \lb{2.108} 
\end{align}
Thus, $M^{H_1}_{1,1}(\lambda^2+ i0,k)=g^{H_1}(\lambda^2+i0,k)$ is
skew-adjoint and hence $H_1$ is reflectionless by \eqref{2.13}.

Similarly, applying \eqref{2.83} yields
\begin{align}
&z\big[M^{H_2}_-(z^2,k)-M^{H_2}_+(z^2,k)\big]^{-1}  \no \\
&\quad =(\chi(k)^{-1})^*\rho(k)^{1/2} 
\big[M^D_-(z,k)+z\rho(k)^{-1}\big]\big[M^{D}_-(z,k)
-M^{D}_+(z,k)\big]^{-1} \no \\ 
&\qquad \times 
\big[M^D_+(z,k)+z\rho(k)^{-1}\big]\rho(k)^{1/2}\chi(k)^{-1}, \quad
k\in\bbZ  \lb{2.109} 
\end{align}
and hence, for all $k\in\bbZ$ and a.e.\
$\lambda\in\sigma_{\ess}(D)$,
\begin{align}
&\lambda\big[M^{H_2}_-(\lambda^2+i0,k)
-M^{H_2}_+(\lambda^2+i0,k)\big]^{-1}  \no \\ 
&\quad =(\chi(k)^{-1})^*\rho(k)^{1/2} 
\big[M^D_-(\lambda+i0,k)+\lambda\rho(k)^{-1}\big] \no \\ 
&\qquad \times
\big[M^{D}_-(\lambda+i0,k)-M^{D}_+(\lambda+i0,k)\big]^{-1} 
\no \\
& \qquad \times \big[M^D_+(\lambda+i0,k)+\lambda\rho(k)^{-1}\big]
\rho(k)^{1/2}\chi(k)^{-1}. 
\lb{2.110} 
\end{align}
Next, consider the Herglotz matrix
\begin{equation}
M_\alpha(z,k)=\big[M^D_-(z,k)+\alpha(k)\big]\big[M^D_-(z,k)
-M^D_+(z,k)\big]^{-1}
\big[M^D_+(z,k)+\alpha(k)\big] \lb{2.111}
\end{equation}
for some self-adjoint $m\times m$ matrix $\alpha(k)$. We claim that 
for all $k\in\bbZ$ and a.e.\ $\lambda\in\sigma_{\ess}(D)$,
\begin{equation}
M_\alpha(\lambda+i0,k)=iC_\alpha(\lambda,k)\, 
\text{ with $C_\alpha(\lambda,k)=C_\alpha(\lambda,k)^*$.} \lb{2.112}
\end{equation}
Indeed, one computes
\begin{align}
M_\alpha(\lambda+i0,k)&=M_-^D(\lambda+i0,k)
\big[M_-^D(\lambda+i0,k)-M_+^D(\lambda+i0,k)\big]^{-1}
M_+^D(\lambda+i0,k) \no \\ 
& \quad +\alpha(k)
\big[M_-^D(\lambda+i0,k)-M_+^D(\lambda+i0,k)\big]^{-1}
\alpha(k) \no \\
& \quad +\alpha(k)
\big[M_-^D(\lambda+i0,k)-M_+^D(\lambda+i0,k)\big]^{-1} 
M_+^D(\lambda+i0,k) \no \\ 
& \quad +M_-^D(\lambda+i0,k) 
\big[M_-^D(\lambda+i0,k)-M_+^D(\lambda+i0,k)\big]^{-1}\alpha(k).
\lb{2.113}
\end{align}
The first two terms on the right-hand side of \eqref{2.113} are
already of the required form \eqref{2.112}. The last two terms on
the right-hand side of \eqref{2.113} can be rewritten in the form
\begin{align}
&\alpha(k)
\big[M_-^D(\lambda+i0,k)-M_+^D(\lambda+i0,k)\big]^{-1} 
M_+^D(\lambda+i0,k) \no \\ 
& \quad +M_-^D(\lambda+i0,k) 
\big[M_-^D(\lambda+i0,k)-M_+^D(\lambda+i0,k)\big]^{-1}
\alpha(k) \no \\
& = (1/2)\alpha(k)\big[M_-^D(\lambda+i0,k)
-M_+^D(\lambda+i0,k)\big]^{-1} \no \\
& \quad \times\big[M_-^D(\lambda+i0,k)+M_+^D(\lambda+i0,k)\big] \no
\\ 
& \quad + (1/2)\big[M_-^D(\lambda+i0,k)+M_+^D(\lambda+i0,k)\big] 
\no \\
& \quad \times \big[M_-^D(\lambda+i0,k)
-M_+^D(\lambda+i0,k)\big]^{-1}\alpha(k) \no \\
&= (1/2)i[\alpha(k)C_{1,2}(\lambda,k)+C_{1,2}(\lambda,k)^*\alpha(k)]
\lb{2.114}
\end{align}
using \eqref{2.105}. Thus, also the last two terms on the
right-hand side of \eqref{2.113} are of the required form
\eqref{2.112} which completes the proof of the claim \eqref{2.112}.
The result \eqref{2.112} applied to \eqref{2.110} then proves that
$M^{H_2}_{1,1}(\lambda^2+ i0,k)=g^{H_2}(\lambda^2+i0,k)$ is
skew-adjoint and hence also $H_2$ is reflectionless.
\end{proof}

The next result is a Borg-type theorem for supersymmetric
Dirac difference operators $D$. However, unlike the Borg-type
theorem for (matrix-valued) Schr\"odinger, Dirac, and Jacobi
operators, this analog for Dirac difference operators
displays a characteristic nonuniqueness feature. 

\begin{theorem}\lb{t4.16}
Assume Hypothesis\ \ref{h4.11} and suppose that $\rho$ and $\chi$ are
reflectionless. Let $D(\rho,\chi)= \left(\begin{smallmatrix} 0 &
\rho S^+ +\chi^* \\ 
\rho^- S^- + \chi& 0 \end{smallmatrix}\right)$ be the associated
self-adjoint Dirac difference operator in 
$\ell^2(\bbZ)^m\oplus\ell^2(\bbZ)^m$
$($cf.~\eqref{2.66}$)$ and suppose that $\sigma(D(\rho,\chi))
=\big[-E_+^{1/2},-E_-^{1/2}\big]\cup \big[E_-^{1/2},E_+^{1/2}\big]$
for some $0\leq E_-<E_+$. Then $\rho$ and $\chi$ are of the form,
\begin{align}
\rho(k)&=\diag(\rho_1(k),\dots,\rho_m(k)), \no \\
\rho_j(k)&=\f{1}{2}\big(E_+^{1/2}-\varepsilon_j E_-^{1/2}\big),
\quad  1\leq j \leq m, \; k\in\bbZ, \lb{2.115} \\ 
\chi(k)&=\diag(\chi_1(k),\dots,\chi_m(k)), \no \\
\chi_j(k)&=\f{1}{2}\big(E_+^{1/2}+\varepsilon_j E_-^{1/2}\big),
\quad  1\leq j \leq m, \; k\in\bbZ, \lb{2.116}
\end{align}
where
\begin{equation}
\varepsilon_j \in \{1,-1\}, \quad 1\leq j\leq m. \lb{2.117}
\end{equation}
\end{theorem}
\begin{proof}
Since $D$ is reflectionless by hypothesis, so are $H_1$ and $H_2$
by Lemma\ \ref{l4.15}. By \eqref{2.54}, 
\begin{equation}
\sigma(D)
=\big[-E_+^{1/2},-E_-^{1/2}\big]\cup \big[E_-^{1/2},E_+^{1/2}\big]
\lb{2.118}
\end{equation}
implies
\begin{equation}
\sigma(H_\ell)=\big[E_-,E_+\big], \quad \ell=1,2. \lb{2.119}
\end{equation}
Applying Theorem\ \ref{t4.7} to $H_1$ and $H_2$ then yields
\begin{align}
\f{1}{4}(E_+-E_-)I_m&=A_1=\rho\chi^+=A_2=\chi\rho, \lb{2.120} \\
\f{1}{4}(E_++E_-)I_m&=B_1=\rho^2+\chi^*\chi=B_2=
(\rho^-)^2+\chi\chi^*, \lb{2.121}
\end{align}
By \eqref{2.120}, $\chi$ satisfies 
\begin{equation}
\chi=\chi^+=\f{1}{4}(E_+-E_-)\rho^{-1} \lb{2.122}
\end{equation}
and hence $\chi$ is a constant (i.e., $k$-independent) positive
definite sequence of $m\times m$ matrices. By \eqref{2.121}
(equivalently, by \eqref{2.122}), then also $\rho=\rho^-$ is a
constant sequence of $m\times m$ matrices. Insertion of
\eqref{2.122} into \eqref{2.121} then yields
\begin{equation}
\rho^4-\f{1}{4}(E_++E_-)\rho^2+\f{1}{16}(E_+-E_-)^2 I_m=0. 
\lb{2.123}
\end{equation}
Since by hypothesis $\rho$ is a positive definite diagonal matrix,
solving the quadratic equation \eqref{2.123} for $\rho_j$, 
$1\leq j\leq m$, yields \eqref{2.115} and hence also \eqref{2.116}
using \eqref{2.122}.
\end{proof}

To the best of our knowledge this result is new even in the
scalar case $m=1$. In particular, the sign ambiguities displayed in
\eqref{2.115} and \eqref{2.116}, giving rise to $2^m$ isospectral
supersymmetric Dirac difference operators, are new in this Borg-type
context. The sign ambiguity disappears and hence uniqueness of the
corresponding inverse spectral problem is restored only in the special
case $E_-=0$, that is, whenever the spectral gap
$\big(-E_-^{1/2},E_-^{1/2}\big)$ of $D(\rho,\chi)$ closes.

Using the supersymmetric formalism described in this section, the
proof of following result can be reduced to that of Theorem\
\ref{t4.9}.

\begin{theorem} \lb{t4.17} 
In addition to assuming Hypothesis\ \ref{h4.11}, suppose that $\rho$
and $\chi$ are periodic with the same period. Let $D(\rho,\chi)=
\left(\begin{smallmatrix} 0 & \rho S^+ +\chi^* \\ 
\rho^- S^- + \chi& 0 \end{smallmatrix}\right)$
be the associated self-adjoint Dirac difference operator in
$\ell^2(\bbZ)^m\oplus\ell^2(\bbZ)^m$ and suppose that
$D(\rho,\chi)$ has uniform spectral multiplicity $2m$. Then
$D(\rho,\chi)$ is reflectionless and hence for all $k\in\bbZ$,
\begin{equation}
\Upsilon^D(\lambda,k)=\f{1}{2} I_{2m} \, \text{ for a.e.\
$\lambda\in\sigma_{\ess}(D(\rho,\chi))$.}  \lb{2.124}
\end{equation}
In particular, assume that $\rho$ and $\chi$ are periodic with the
same period, that $D(\rho,\chi)$ has uniform spectral multiplicity
$2m$, and that $\sigma(D(\rho,\chi))=[-E_+^{1/2},-E_-^{1/2}] \cup 
[E_-^{1/2},E_+^{1/2}]$ for some $0\leq E_-<E_+$. Then
$\rho$ and $\chi$ are of the form,
\begin{align}
\rho(k)&=\diag(\rho_1(k),\dots,\rho_m(k)), \no \\
\rho_j(k)&=\f{1}{2}\big(E_+^{1/2}-\varepsilon_j E_-^{1/2}\big),
\quad  1\leq j \leq m, \; k\in\bbZ, \lb{2.125} \\ 
\chi(k)&=\diag(\chi_1(k),\dots,\chi_m(k)), \no \\
\chi_j(k)&=\f{1}{2}\big(E_+^{1/2}+\varepsilon_j E_-^{1/2}\big),
\quad  1\leq j \leq m, \; k\in\bbZ, \lb{2.126}
\end{align}
where
\begin{equation}
\varepsilon_j \in \{1,-1\}, \quad 1\leq j\leq m. \lb{2.127}
\end{equation}
\end{theorem}



\begin{thebibliography}{99}
%
\bi{AG94} D.~Alpay and I.~Gohberg, {\it Inverse spectral problems
for difference operators with rational scattering matrix function},
Integr. Equat. Oper. Th. {\bf 20}, 125--170 (1994).
%
\bi{AN84} A.~I.~Aptekarev and E.~M.~Nikishin, {\it The scattering
problem for a discrete Sturm-Liouville problem}, Math. USSR
Sborn. {\bf 49}, 325--355 (1984).
%
\bi{AD56} N.~Aronszajn and W.~F.~Donoghue, {\it On exponential
representations of analytic functions in the upper half-plane with 
positive imaginary part}, J. Anal. Math. {\bf 5}, 321-388 (1956-57). 
%
\bi{Be68} Ju.~Berezanskii, {\it Expansions in Eigenfunctions
of Selfadjoint Operators}, Transl. Math. Mongraphs, Vol.\ 17,
Amer. Math. Soc., Providence, R.I., 1968.
%
\bi{BGS86} Yu.~M.~Berezanskii, M.~I.~Gekhtman, and
M.~E.~Shmoish, {\it Integration of some chains of nonlinear
difference equations by the method of the inverse spectral
problem}, Ukrain. Math. J. {\bf 38}, 74--78 (1986).
%
%
\bi{BG90} Yu.~M.~Berezanskii and M.~I.~Gekhtman, {Inverse problem
of the spectral analysis and non-abelian chains of nonlinear
equations}, Ukrain. Math. J. {\bf 42}, 645--658 (1990).
%
\bi{Bo46} G.~Borg, {\it Eine Umkehrung der Sturm-Liouvilleschen
Eigenwertaufgabe}, Acta Math. {\bf 78}, 1--96 (1946). 
%
\bi{BGHT98}W.~Bulla, F.~Gesztsy, H.~Holden, and G.~Teschl, {\it Algebro-
Geometric Quasi-Periodic Finte-Gap Solutions of the Toda and Kac-van
Moerbeke Hierarachies,} Memoirs of the Amer. Math. Soc. {\bf 135/641},
(1998).
%
\bi{Ca76} R.~W.~Carey, {\it A unitary invariant for pairs of self-adjoint
operators}, J. reine angew. Math. {\bf 283}, 294--312 (1976).
%
\bi{CG01} S.~Clark and F.~Gesztesy, {\it Weyl--Titchmarsh
$M$-function asymptotics for matrix-valued Schr\"odinger
operators}, Proc. London Math. Soc. {\bf 82}, 701--724 (2001).
%
\bi{CG02} S.~Clark and F.~Gesztesy, {\it Weyl--Titchmarsh
$M$-function asymptotics and Borg-type theorems for Dirac
operators}, Trans. Amer. Math. Soc. {\bf 354}, 3475--3534 (2002).
%
\bi{CG03} S.~Clark and F.~Gesztesy, {\it On Weyl--Titchmarsh Theory 
for Singular Finite Difference Hamiltonian Systems}, J. Comput. Appl.
Math., to appear.
%
\bi{CGHL00} S.~Clark, F.~Gesztesy, H.~Holden, and
B.~M.~Levitan, {\it Borg-type theorems for matrix-valued
Schr\"odinger operators}, J. Diff. Eqs. {\bf 167}, 181--210 (2000).
%
\bi{De78} P.\ A.\ Deift, {\it Applications of a commutation formula},
Duke Math. J. {\bf 45}, 267--310 (1978).
%
\bi{De95} B.~Despr\'es, {\it The Borg theorem for the vectorial Hill's
equation}, Inverse Probl. {\bf 11}, 97--121 (1995). 
%
\bi{DL96} A.~J.~Duran and P.~Lopez-Rodriguez, {\it Orthogonal
matrix polynomials: zeros and Blumenthal's theorem}, J. Approx. Th.
{\bf 84}, 96--118 (1996).
%
\bi{DL00} A.~J.~Duran and P.~Lopez-Rodriguez, {\it $N$-extremal
matrices of measures for an indeterminate matrix moment problem}, J.
Funct. Anal. {\bf 174}, 301--321 (2000).
%
\bi{DV95} A.~J.~Duran and W.~ Van Assche, {\it Orthogonal matrix
polynomials and higher-order recurrence relations}, Lin. Algebra
Appl. {\bf 219}, 261--280 (1995).
%
\bi{Fl75} H.\ Flaschka, {\it Discrete and periodic illustrations of
some aspects of the inverse method}, in {\it Dynamical Systems, Theory
and Applications}, J.\ Moser (ed.), Lecture Notes In Physics, Vol.\ 38,
Springer Verlag, Berlin, 1975, p.\ 441--466.
%
\bi{Fu76} M.~Fukushima, {\it A spectral representation on ordinary
linear difference equation with operator-valued coefficients of the
second order}, J. Math. Phys. {\bf 17}, 1084--1072 (1976).
%
\bi{Ge82} J.~S.~Geronimo, {\it Scattering theory and matrix
orthogonal polynomials on the real line}, Circuits Syst. Signal
Process. {\bf 1}, 471--495 (1982).
%
\bi{GH05} F.\ Gesztesy and H.\ Holden, {\it Soliton Equations and
Their Algebro-Geometric Solutions. Vol. II: $(1+1)$-Dimensional 
Discrete Models},  Cambridge Studies in Advanced Mathematics,
Cambridge Univ. Press, in preparation.
%
\bi{GKM02} F.~Gesztesy, A.~Kiselev, and K.~A.~Makarov,
{\it Uniqueness Results for Matrix-Valued Schr\"odinger,
Jacobi, and Dirac-Type Operators}, Math. Nachr. {\bf 239--240},
103--145 (2002).
%
\bi{GKT96} F.~Gesztesy, M.~Krishna, and G.~Teschl, {\it On
isospectral sets of Jacobi operators,} Commun. Math. Phys.
{\bf 181}, 631--645 (1996).
%
\bi{GSS91} F.\ Gesztesy, W.\ Schweiger, and B.\ Simon, {\it Commutation
methods applied to the mKdV-equation}, Trans. Amer. Math. Soc. {\bf
324}, 465--525 (1991).
%
\bi{GS96} F.~Gesztesy and B.~Simon, {\it The $\xi$ function},
Acta Math. {\bf 176}, 49--71 (1996).
%
\bibitem{GS96a} F.~Gesztesy and B.~Simon, {\it Uniqueness
theorems in inverse spectral theory for one-dimensinal Schr\"odinger
operators}, Trans. Amer. Math. Soc. {\bf 348}, 349--373 (1996).
%
\bi{GS97} F.~Gesztesy and B.~Simon, {\it $m$-functions and inverse 
spectral analysis for finite and semi-infinite Jacobi matrices},
J. Analyse Math. {\bf 73}, 267--297 (1997).
%
\bi{GS99} F.~Gesztesy and B.~Simon, {\it On local
Borg-Marchenko uniqueness results}, Commun. Math. Phys.
{\bf 211}, 273--287 (2000).
%
\bi{GT00} F.~Gesztesy and E.~Tsekanovskii, {\it On
matrix-valued Herglotz functions}, Math. Nachr. {\bf 218}, 61--138
(2000).
%
\bi{HS81} D.~B.~Hinton and J.~K.~Shaw, {\it On Titchmarsh--Weyl
$M(\lambda)$-functions for linear Hamiltonian systems},
J. Diff. Eqs. {\bf 40}, 316--342 (1981).
%
\bi{HS82} D.~B.~Hinton and J.~K.~Shaw, {\it On the spectrum
of a singular Hamiltonian system}, Quaest. Math.
{\bf 5}, 29--81 (1982).
%
\bi{HS86} D.~B.~Hinton and J.~K.~Shaw, {\it On the
spectrum of a singular Hamiltonian system, II}, Quaest. Math. {\bf 10},
1--48 (1986).
%
\bi{KS03} R.\ Killip and B.\ Simon, {\it Sum rules for Jacobi matrices
and their applications to spectral theory}, Ann. Math. {\bf 158},
253--321 (2003).
%
\bi{KM98} A.~G.~Kostyuchenko and K.~A.~Mirzoev, {\it Three-term
recurrence relations with matrix coefficients. The completely indefinite
case}, Math. Notes {\bf 63}, 624--630 (1998).
%
\bi{KM99} A.~G.~Kostyuchenko and K.~A.~Mirzoev, {\it Generalized Jacobi
matrices and deficiency numbers of ordinary differential
operators with polynomial coefficients}, Funct. Anal. Appl. {\bf 33},
25--37 (1999).
%
\bi{KM01} A.~G.~Kostyuchenko and K.~A.~Mirzoev, {\it Complete
indefiniteness tests for Jacobi matrices with matrix entries}, Funct.
Anal. Appl. {\bf 35}, 265--269 (2001).
%
\bi{KS88} S.~Kotani and B.~Simon, {\it Stochastic
Schr\"odinger
operators and Jacobi matrices on the strip}, Commun. Math.
Phys. {\bf 119}, 403--429 (1988).
%
\bi{Lo99} P.~L{\'o}pez-Rodriguez, {\it Riesz's theorem for
orthogonal matrix polynomials}, Constr. Approx. {\bf 15},
135--151 (1999).
%
\bi{MBO92} F.~G.~Maksudov, E.~M.~Bairamov, and R.~U.~Orudzheva, {\it
The inverse scattering problem for an infinite Jacobi matrix with
operator elements}, Russ. Acad. Sci. Dokl. Math. {\bf 45},
366--370 (1992).
%
\bi{Ma94} M.~M.~Malamud, {\it Similarity of Volterra 
operators and
related questions of the theory of differential equations 
of fractional
order}, Trans.~Moscow Math.~Soc. {\bf 55}, 57--122 (1994). 
%
\bi{Ma99} M.~M.~Malamud,
{\it Uniqueness questions in inverse problems for systems of
ordinary
differential equations on a finite interval}, Trans.~Moscow
Math.~Soc. {\bf 60}, 173--224 (1999).
%
\bi{Ma99a} M.~M.~Malamud,
{\it Borg type theorems for first-order systems on a finite 
interval}, Funct. Anal. Appl. {\bf 33}, 64--68 (1999).
%
\bi{Os97} A.~S.~Osipov, {\it Integration of non-abelian Lanmuir type
lattices by the inverse spectral problem method}, Funct. Anal. Appl. {\bf
31}, 67--70 (1997). 
%
\bi{Os00} A.~S.~Osipov, {\it Some properties of resolvent sets of
second-order difference operators with matrix coefficients}, Math. Notes
{\bf 68}, 806--809 (2000). 
%
\bi{Os02} A.~Osipov, {\it On some issues related to the moment problem
for the band matrices with operator elements}, J. Math. Anal. Appl. {\bf
275}, 657--675 (2002).
%
\bi{RS02} J.\ Rovnyak and L.\ A.\ Sakhnovich, {\it Some indefinite 
cases of spectral problems for canonical systems of difference
equations}, Lin. Algebra Appls. {\bf 343--344}, 267--289 (2002).
%
\bi{Sa03} A.\ L.\ Sakhnovich, {\it Discete canonical systems and
non-Abelian Toda lattice: B\"acklund--Darboux transformation, 
Weyl-functions, and explicit solutions}, preprint, 2004.
%
\bi{Sa97}  L.~A.~Sakhnovich, {\it Interpolation Theory and its
Applications}, Kluwer, Dordrecht, 1997.
%
\bi{Se80} V.~P.~Serebrjakov, {\it The inverse problem of scattering
theory for difference equations with matrix coefficients}, Sov. Math.
Dokl. {\bf 21}, 148--151 (1980).
%
\bi{SY95} M.~Sodin and P.~Yuditskii, {\it Almost periodic Sturm-Liouville
operators with Cantor  homogeneous spectrum}, Comment. Math. Helvetici 
{\bf 70}, 639--658 (1995).
%
\bi{SY96} M.~Sodin and P.~Yuditskii, {\it Almost periodic Jacobi
matrices with homogeneous spectrum, infinite dimensional Jacobi
inversion, and Hardy spaces of character-automorphic functions}, J.
Geom. Anal. {\bf 7}, 387--435 (1997).
%
\bi{Te98} G.~Teschl, {\it Trace formulas and inverse
spectral theory for Jacobi operators}, Commun. Math. Phys.
{\bf 196}, 175--202 (1998).
%
\bi{Te00} G.~Teschl, {\it Jacobi Operators and Completely Integrable
Nonlinear Lattices}, Math. Surv. Monographs, Vol.\ 72,
Amer. Math. Soc., Providence, R.I., 2000.
%
\end{thebibliography}
\end{document}